\documentclass[10pt]{amsart}
\usepackage[T1]{fontenc}
\usepackage{a4wide}
\usepackage[latin1]{inputenc}
\usepackage{amsmath}
\usepackage{amsfonts}
\usepackage{amssymb}
\usepackage{amsthm}
\usepackage{subfigure}
\newtheorem{theo}{Theorem}[section]
\newtheorem{prop}[theo]{Proposition}
\newtheorem{coro}[theo]{Corollary}
\newtheorem{lemm}[theo]{Lemma}
\theoremstyle{definition}

\theoremstyle{remark}
\newtheorem*{rema}{Remark}
\newcommand{\Op}{\text{Op}}

\title{Entropy of semiclassical measures for nonpositively curved surfaces}

\author[G. Rivi\`ere]{Gabriel Rivi\`ere}
\address{Centre de Math\'ematiques Laurent Schwartz, \'Ecole Polytechnique, 91128 Palaiseau Cedex, France}
\email{gabriel.riviere@math.polytechnique.fr}

\begin{document}

\maketitle

\begin{abstract}
We study the asymptotic properties of eigenfunctions of the Laplacian in the case of a compact Riemannian surface of nonpositive sectional curvature. We show that the Kolmogorov-Sinai entropy of a semiclassical measure $\mu$ for the geodesic flow $g^t$ is bounded from below by half of the Ruelle upper bound, i.e.
$$h_{KS}(\mu,g)\geq \frac{1}{2}\int_{S^*M} \chi^+(\rho) d\mu(\rho).$$
We follow the same main strategy as in~\cite{GR} and refer the reader to it for the details of several lemmas.
\end{abstract}

\section{Introduction}

Let $M$ be a compact $\mathcal{C}^{\infty}$ Riemannian manifold. For all $x\in M$, $T^*_xM$ is endowed with a norm $\|.\|_x$ given by the metric over $M$. The geodesic flow $g^t$ over $T^*M$ is defined as the Hamiltonian flow corresponding to the Hamiltonian $H(x,\xi):=\frac{\|\xi\|_x^2}{2}$. This last quantity corresponds to the classical kinetic energy in the case of the absence of potential. As any observable, this quantity can be quantized via pseudodifferential calculus and the quantum operator corresponding to $H$ is $-\frac{\hbar^2\Delta}{2}$ where $\hbar$ is proportional to the Planck constant and $\Delta$ is the Laplace Beltrami operator acting on $L^2(M)$. Our main concern in this note will be to study the asymptotic behavior, as $\hbar$ tends to $0$, of the following sequence of distributions:
$$\forall a\in\mathcal{C}^{\infty}_o(T^*M),\ \mu_{\hbar}(a)=\int_{T^*M}a(x,\xi)d\mu_{\hbar}(x,\xi):=\langle \psi_{\hbar},\Op_{\hbar}(a)\psi_{\hbar}\rangle_{L^2(M)},$$
where $\Op_{\hbar}(a)$ is a $\hbar$-pseudodifferential operator of symbol $a$~\cite{DS} and $\psi_{\hbar}$ satisfies
$$-\hbar^2\Delta\psi_{\hbar}=\psi_{\hbar}.$$
An accumulation point of such a sequence of distribution $\mu_{\hbar}$ is called a semiclassical measure. Moreover, one knows that a semiclassical measure is a probability measure on $S^*M:=\{\|\xi\|_x^2=1\}$ which is invariant under the geodesic flow $g^t$ on $S^*M$. For manifolds of negative curvature, the geodesic flow on $S^*M$ satisfies strong chaotic properties (Anosov property, ergodicity of the Liouville measure) and as a consequence, it can be shown that \emph{almost all} the $\mu_{\hbar}$ converge to the Liouville measure on $S^*M$~\cite{Sc},~\cite{Ze},~\cite{CdV}. This phenomenon is known as the quantum ergodicity property. A main challenge concerning this result would be to answer the Quantum Unique Ergodicity Conjecture~\cite{RS}, i.e. determine whether the Liouville measure is the only semiclassical measure or not (at least for manifolds of negative curvature).\\
In~\cite{An}, Anantharaman used the Kolmogorov-Sinai entropy to study the properties of semiclassical measures on manifolds of negative curvature\footnote{In fact, her result was about manifolds with Anosov geodesic flow, for instance manifolds of negative curvature.}. In particular, she showed that the Kolmogorov-Sinai entropy of any semiclassical measure is positive. This result implies that the support of a semiclassical measure cannot be restricted to a closed geodesic, i.e. \emph{eigenfunctions of the Laplacian cannot concentrate only on a closed geodesic in the high energy limit}. In subsequent works, with Nonnenmacher and Koch, more quantitative lower bounds on the entropy of semiclassical measures were given~\cite{AN2},~\cite{AKN}.

\subsection{Kolmogorov-Sinai entropy}

Let us recall a few facts about the Kolmogorov-Sinai (also called metric) entropy (see~\cite{Wa} or appendix~\ref{KSentropy} for more details and definitions). It is a nonnegative number associated to a flow $g$ and a $g$-invariant measure $\mu$, that estimates the complexity of $\mu$ with respect to this flow. For example, a measure carried by a closed geodesic will have entropy zero. Recall also that a standard theorem of dynamical systems due to Ruelle~\cite{R} asserts that, for any invariant measure $\mu$ under the geodesic flow:
\begin{equation}\label{ruelle}h_{KS}(\mu,g)\leq \int_{S^*M} \sum_{j}\chi^+_j(\rho) d\mu(\rho)\end{equation}
with equality if and only if $\mu$ is the Liouville measure in the case of an Anosov flow~\cite{LY}. In the previous inequality, the $\chi_j^+$ denoted the positive Lyapunov exponents of $(S^*M,g^t,\mu)$~\cite{BP}.\\
Regarding these properties, the main result of Anantharaman-Koch-Nonnenmacher was to show that, for a semiclassical measure $\mu$ on an Anosov manifold, one has
$$h_{KS}(\mu,g)\geq \int_{S^*M} \sum_{j=1}^{d-1}\chi^+_j(\rho) d\mu(\rho)-\frac{(d-1)\lambda_{\max}}{2}.$$
where $\lambda_{\max}:=\lim_{t\rightarrow\pm\infty}\frac{1}{t}\log\sup_{\rho\in S^*M}|d_{\rho}g^t|$ is the maximal expansion rate of the geodesic flow and the $\chi^+_j$'s are the positive Lyapunov exponents~\cite{BP}. Compared with the original result from~\cite{An}, this inequality gives a precise lower bound on the entropy of a semiclassical measure. For instance, for manifolds of \emph{constant} negative curvature, this lower bound can be rewritten as $\frac{d-1}{2}$. However, it can turn out that $\lambda_{\max}$ is a very large quantity and in this case, the previous lower bound can be negative (which would imply that it is an empty result). Combining these two observations~\cite{AN2}, they were lead to formulate the conjecture that, for any semiclassical measure $\mu$, one has
$$h_{KS}(\mu,g)\geq \frac{1}{2}\int_{S^*M} \sum_{j=1}^{d-1}\chi^+_j(\rho) d\mu(\rho).$$
They also ask about the extension of this conjecture to manifolds without conjugate points~\cite{AN2}. In recent work~\cite{GR}, we were able to prove that their conjecture holds for any surface with an Anosov geodesic flow (for instance surfaces of negative curvature). Regarding our proof and the nice properties of surfaces of nonpositive curvature~\cite{Ru},~\cite{Eb}, it became clear that our result can be adapted in the following way:
\begin{theo}\label{maintheo2} Let $M$ be a $\mathcal{C}^{\infty}$ Riemannian surface of nonpositive sectional curvature and $\mu$ a semiclassical measure. Then,
\begin{equation}\label{mainineq2}h_{KS}(\mu,g)\geq \frac{1}{2}\int_{S^*M} \chi^+(\rho) d\mu(\rho),\end{equation}
where $h_{KS}(\mu,g)$ is the Kolmogorov-Sinai entropy and $\chi^+(\rho)$ is the upper Lyapunov exponent at point $\rho$.
\end{theo}
In particular, this result shows that the support of a semiclassical measure cannot be reduced to a closed unstable geodesic. We underline that our inequality is also coherent with the quasimodes constructed by Donnelly~\cite{Do}. In fact, his quasimodes are supported on closed stable geodesics (included in flat parts of a surface of nonpositive curvature) and have zero entropy. We can make a last observation on the assumptions on the manifold: it is not known whether the geodesic flow is ergodic or not for the Liouville measure on a surface of nonpositive curvature. The best result in this direction is that there exists an open invariant subset $U$ of positive Liouville measure such that the restriction $g_{|U}$ is ergodic with respect to Liouville~\cite{BP}. The extension of this result on the entropy of semiclassical measures raises the question of knowing whether one could obtain an analogue of this result for weakly chaotic systems. For instance, regarding the counterexamples constructed in~\cite{Has}, it would be interesting to have a lower bound for ergodic billiards.\\
Our purpose in this note is to prove theorem~\ref{maintheo2}. Our strategy will be the same as in~\cite{GR}. So we will focus on the main differences and refer the reader to~\cite{GR} and~\cite{AN2} for the details of several lemmas. The crucial observation is that as in the Anosov case, surfaces of nonpositive curvature have \emph{continuous stable and unstable foliations} and \emph{no conjugate points}. This property was at the heart of the proofs in~\cite{AN2},~\cite{AKN} and~\cite{GR} and we will verify that even if the properties of these stable/unstable directions are weaker for surfaces of nonpositive curvature, they are sufficient to prove the conjecture of Anantharaman-Nonnenmacher in this weakly chaotic setting. In~\cite{AN2},~\cite{AKN},~\cite{GR}, there was a dynamical quantity which was crucially used: the unstable Jacobian of the geodesic flow. In the case of surfaces of nonpositive curvature, one can introduce an analogue of it. This quantity comes from the study of Jacobi fields and is called the unstable Riccati solution $U^u(\rho)$~\cite{Ru}. In this setting, it has been shown that the Ruelle inequality can be rewritten as follows~\cite{FM}:
$$h_{KS}(\mu,g)\leq \int_{S^*M} U^u(\rho) d\mu(\rho).$$
So, the lower bound of theorem~\ref{maintheo2} can be rewritten as
\begin{equation}\label{mainineq3}h_{KS}(\mu,g)\geq \frac{1}{2}\int_{S^*M} U^u(\rho) d\mu(\rho).\end{equation}
The main adavantage of this new formulation is that the function in the integral of the lower bound is defined everywhere (and not almost everywhere).
\begin{rema} One could also ask whether it would be possible to extend this result to surfaces without conjugate points. In fact, these surfaces also have a stable and unstable foliations (and of course no conjugate points). Moreover, according to Green~\cite{Gr} and Eberlein~\cite{Eb0}, the Jacobi fields also satisfy a property of uniform divergence (at least in dimension $2$). The main difficulty is that the continuity of $U^u(\rho)$ is not true anymore~\cite{BBB} and at this point, we do not see any way of escaping this difficulty.
\end{rema}

\subsection{Organization of the article}

In section~\ref{s:classical-setting}, we will give a precise survey\footnote{We refer the reader to~\cite{Eb} or~\cite{Ru} for more details.} on surfaces of nonpositive curvature and highlight the properties we will need to make the proof work. Then, in section~\ref{s:outline}, we will draw a precise outline of the proof and we will refer to~\cite{GR} for the details of some lemmas. In section~\ref{s:main-est}, we will explain how the main result from~\cite{AN2} can be adapted in the setting of surfaces of nonpositive curvature. In section~\ref{s:appl-uncert}, we follow the same strategy as in~\cite{GR} to derive a crucial estimate on the quantum pressures. Finally, in the appendix, we recall some results on quantum pressure from~\cite{AN2}.

\subsection*{Acknowledgements} I would like to sincerely thank my advisor Nalini Anantharaman for introducing me to this question and for encouraging me to extend the result from~\cite{GR} to nonpositively curved surfaces. I also thank her for many helpful discussions about this subject.

\section{Classical setting of the article}\label{s:classical-setting}

\subsection{Surfaces of nonpositive curvature}\label{nonpositive} In this first section, we recall some facts about nonpositively curved manifolds~\cite{Ru},~\cite{Eb}.

\subsubsection{Stable and unstable Jacobi fields} We define $\pi:S^*M\rightarrow M$ the canonical projection $\pi(x,\xi):=x$. The vertical subspace $V_{\rho}$ at the point $\rho=(x,\xi)$ is the kernel of the application $d_{\rho}\pi$. We underline that it is in fact the tangent space in $\rho$ of the $1$-dimensional submanifold $S_x^*M$. In the case of a surface, it has dimension $1$. We can also define the horizontal subspace in $\rho$. Precisely, for $Z\in T_{\rho}S^*M$, we consider a smooth curve $c(t)=(a(t),b(t))$, $t\in(-\epsilon,\epsilon)$, in $S^*M$ such that $c(0)=\rho$ and $c'(0)=Z$. Then, we define the horizontal space $\mathcal{H}_{\rho}$ as the kernel of the application $\mathbb{K}_{\rho}(Z)=\nabla_{a'(0)}b(0)=\nabla_{d_{\rho}\pi(Z)}b(0)$, where $\nabla$ is the Levi-Civita connection. This subspace contains $X_{H}(\rho)$ the vector field tangent to the Hamiltonian flow. For a surface, this subspace is of dimension $2$. We know that we can use these two subspaces to split the tangent space $T_{\rho}S^*M=\mathcal{H}_{\rho}\oplus V_{\rho}$ (it is the usual way to split the tangent space in order to define the Sasaki metric on $S^*M$~\cite{Ru}). Using this decomposition, we would like to recall an important link between the linearization of the geodesic flow and the Jacobi fields on $M$. To do this, we underline that to each point $\rho$ in $S^*M$ corresponds a unique unit speed geodesic $\gamma_{\rho}$. Then we define a Jacobi field in $\rho$ (or along $\gamma_{\rho}$) as a solution of the differential equation:
$$\mathbb{J}"(t)+R(\gamma_{\rho}'(t),\mathbb{J}(t))\gamma_{\rho}'(t)=0,$$
where $R(X,Y)Z$ is the curvature tensor applied to the vector fields $X$, $Y$ and $Z$ and $\mathbb{J}'(t)=\nabla_{\gamma_{\rho}'(t)}\mathbb{J}(t)$. We recall that we can interpret Jacobi fields as geodesic variation vector fields~\cite{Eb}. Precisely, consider a $\mathcal{C}^{\infty}$ family of curves $c_s:[a,b]\rightarrow M$, $s$ in $(-\epsilon,\epsilon)$. We say that it is a $\mathcal{C}^{\infty}$ variation of $c=c_0$. It defines a corresponding variation vector field $Y(t)=\frac{\partial}{\partial s}(c_s(t))_{|s=0}$ that gives the initial velocity of $s\mapsto c_s(t)$. If we suppose now that $c$ is a geodesic of $M$, then a $\mathcal{C}^2$ vector field $Y(t)$ on $c$ is a Jacobi vector field if and only if $Y(t)$ is the the variation vector field of a geodesic variation of $c$ (i.e. $\forall s\in(-\epsilon,\epsilon)$, $c_s$ is a geodesic of $M$). For instance, $\gamma_{\rho}'(t)$ and $t\gamma_{\rho}'(t)$ are Jacobi vector fields along $\gamma_{\rho}$.\\
Consider now a vector $(V,W)$ in $T_{\rho}S^*M$ given in the coordinates $\mathcal{H}_{\rho}\oplus V_{\rho}$. Using the canonical identification given by $d_{\rho}\pi$ and $\mathbb{K}_{\rho}$, there exists a unique Jacobi field $\mathbb{J}_{V,W}(t)$ in $\rho$ whose initial conditions are $\mathbb{J}_{V,W}(0)=V$ and $\mathbb{J}_{V,W}'(0)=W$, such that $$d_{\rho}g^t(V,W)=(\mathbb{J}_{V,W}(t),\mathbb{J}_{V,W}'(t))$$ in coordinates $\mathcal{H}_{g^t\rho}\oplus V_{g^t\rho}$~\cite{Ru} (lemma $1.4$). We define $N_{\rho}$ the subspace of $T_{\rho}S^*M$ of vectors orthogonal to $X_H(\rho)$ and $H_{\rho}$ the intersection of this subspace with $\mathcal{H}_{\rho}$. Using the previous property about Jacobi fields, we know that the subbundle $\mathcal{N}$ perpendicular to the Hamiltonian vector field is invariant by $g^t$ and that we have the following splitting~\cite{Ru} (lemma $1.5$):
$$T_{\rho}S^*M=\mathbb{R}X_H(\rho)\oplus H_{\rho}\oplus V_{\rho}.$$
Obviously, these properties can be extended to any energy layer $\mathcal{E}(\lambda)$ for any positive $\lambda$. Following~\cite{Ru} (lemma $3.1$), we can make the following construction of two particular Jacobi fields along $\gamma_{\rho}$. We denote $(\gamma_{\rho}'(t),e(t))$ an orthonormal basis defined along $\gamma_{\rho}(t)$. Given a positive $T$ and because there are no conjugate points on the manifold $M$, there exists a unique Jacobi field $\mathbb{J}_T(t)$ such that $\mathbb{J}_T(0)=e(0)$ and $\mathbb{J}_T(T)=0$. Moreover, $\mathbb{J}_T(t)$ is perpendicular to $\gamma_{\rho}(t)$ for all $t$ in $\mathbb{R}$~\cite{Ru} (page $50$). As a consequence, $\mathbb{J}_T(t)$ can be identified with its coordinate along $e(t)$ (as $T_{\gamma_{\rho}(t)}M$ is of dimension $2$). A result due to Hopf (lemma $3.1$ in~\cite{Ru}) tells us that the limits
$$\lim_{T\rightarrow+\infty}\mathbb{J}_T(t)\ \text{and}\ \lim_{T\rightarrow-\infty}\mathbb{J}_T(t)$$
exist. They are denoted $\mathbb{J}^s_{\rho}(t)$ and $\mathbb{J}^u_{\rho}(t)$ (respectively the stable and the unstable Jacobi field). They satisfy the simplified one dimensional Jacobi equation:
$$\mathbb{J}"(t)+K(t)\mathbb{J}(t)=0,$$
where $K(t)=K(\gamma_{\rho}(t))$ is the sectional curvature at $\gamma_{\rho}(t)$. They are never vanishing Jacobi fields with $\mathbb{J}^*_{\rho}(0)=e(0)$ and for all $t$ in $\mathbb{R}$, they are perpendicular to $\gamma_{\rho}'(t)$. Moreover, we have $\|\mathbb{J}^{*'}_{\rho}(t)\|\leq\sqrt{K_0}\|\mathbb{J}^{*}_{\rho}(t)\|$ for every $t$ in $\mathbb{R}$ (where $-K_0$ is some negative lower bound on the curvature). Using the previous link between geodesic flow and Jacobi fields, we can lift these subspaces to invariant subspaces $E^s(\rho)$ and $E^u(\rho)$ called the Green stable and unstable subspaces. These subspaces have dimension $1$ (in the case of surfaces) and are included in $N_{\rho}$. A basis of $E^s(g^t\rho)$ is given by $(\mathbb{J}^s_{\rho}(t),\mathbb{J}^{s'}_{\rho}(t))$ in coordinates $H_{g^t\rho}\oplus V_{g^t\rho}$. We can underline that both subspaces are uniformly transverse to $V_{\rho}$ and that it can happen that they are equal to each other (which was not the case in the Anosov setting). In the case of nonpositive curvature, these subspaces depend continuously in $\rho$ and are integrable as in the Anosov case~\cite{Eb}.

\subsubsection{Riccati equation}
In the case where the Green subspaces attached to $\rho$ are linearly independent, a splitting of $N_{\rho}$ is given by $E^{u}(\rho)\oplus E^s(\rho)$ and the splitting holds for all the trajectory. For the opposite case, we know that the Green subspaces attached to $\rho$ (and hence to a geodesic $\gamma_{\rho}$) are linearly dependent if and only if the sectional curvature is vanishing at every point of the geodesic $\gamma_{\rho}$~\cite{Ru}. As a consequence, we cannot use the same kind of splitting. However, there exists a splitting of $N_{\rho}$ that we can use in both cases, precisely $E^{u}(\rho)\oplus V_{\rho}$. We would like to mention that the one dimensional Jacobi equation defined earlier gives rise to the Riccati equation:
$$U'(t)+U^2(t)+K(t)=0,$$
where $U(t)=\mathbb{J}'(t)\mathbb{J}(t)^{-1}$ for non vanishing $\mathbb{J}$. Then we define the corresponding unstable Riccati solution associated to the unstable Jacobi field as $U^u_{\rho}(t):=\mathbb{J}^{u'}_{\rho}(t)(\mathbb{J}^u_{\rho}(t))^{-1}$. It is a nonnegative quantity and it describes the growth of the unstable Jacobi field (in dimension $2$) as follows:
$$\|\mathbb{J}^{u}_{\rho}(t)\|=\|\mathbb{J}^{u}_{\rho}(0)\|e^{\int_0^tU^u_{\rho}(s)ds}.$$
The same works for the stable Jacobi field. Both quantities are continuous\footnote{The continuity in $\rho$ is a crucial property that we will use in our proof. We underline that it is not true if we only suppose the surface to be without conjugate points~\cite{BBB}.} with respect to $\rho$. We underline that, we can use the previous results to obtain the bound $\|d_{\rho}g^t_{|E^{u}(\rho)}\|\leq \sqrt{1+K_0}e^{\int_0^tU^u_{\rho}(s)ds}$. So the unstable Riccati solution describe the infinitesimal growth of the geodesic flow along the unstable direction, whereas $J^u(\rho)^{-1}$ used in the previous sections described the growth at time $1$. More precisely, as for the unstable Jacobian, Freire and Ma\~{n}\'e showed that the unstable Riccati solutions are related to the Lyapunov exponents. In fact, they proved that the Ruelle bound for the entropy of a $g$-invariant measure $\mu$ in the case of nonpositive curvature (precisely for manifolds without conjugate points)~\cite{FM} is:
$$h_{KS}(\mu,g)\leq \int_{S^*M}U^u(\rho)d\mu(\rho).$$

\subsubsection{Divergence of vanishing Jacobi fields}
A last point we would like to recall is a result due to Green~\cite{Gr} and to Eberlein in the general case~\cite{Eb0}. It asserts that for any positive $c$ there exists a positive $T=T(c)$ such that for any $\rho$ in $S^*M$ and for any nontrivial Jacobi field $\mathbb{J}(t)$ along $\gamma_{\rho}$ such that $\mathbb{J}(0)=0$ and $\|\mathbb{J}'(0)\|\geq 1$, for all $t$ larger than $T$, we have $\|\mathbb{J}(t)\|\geq c$ (proposition $3.1$~\cite{Ru}). This property of uniform divergence only holds in dimension $2$ and as it is crucially used in the following, our proof only works for surfaces of nonpositive curvature. In larger dimensions, the same result holds but without any uniformity in $\rho$. Finally, all these properties allow to prove the following lemma:
\begin{lemm}\label{unstable} Let $v=(0,V)$ be a unit vertical vector at $\rho$. Then for any $c>0$, there exists $T=T(c)>0$ (independent of $\rho$ and of $v$) such that for any $t\geq T$, $\|d_{\rho}g^tv\|\geq c$.
\end{lemm}
We underline that, for $t\geq T$, the angle between $E^{u}(g^t\rho)$ and $d_{\rho}g^tv$ is bounded by some $\kappa(c)$ with $\kappa(c)$ arbitrarly small as $c$ tends to infinity.

\subsection{Discretization of the unstable Riccati solution}\label{classical-nonpos}

For $\theta$ small positive number ($\theta$ will be fixed all along the paper), one defines $$\mathcal{E}^{\theta}:=H^{-1}(]1/2-\theta,1/2+\theta[).$$
From previous section, we know that there exists a constant $b_0$ such that
$$\forall\rho\in \mathcal{E}^{\theta},\ 0\leq U^u(\rho)\leq b_{0}.$$
This function will replace the logarithm of the unstable Jacobian $\log J^u$ in the proof from~\cite{GR}. The situation is slightly different from the case of an Anosov flow as we do not have that $U^u$ is uniformly bounded from below by some positive constant, a property that was crucially to prove theorem $1.2$ in~\cite{GR}. We solve this problem by introducing a small positive parameter $\epsilon_0$ and defining an auxiliary function
$$U_0^u(\rho):=\sup\{U^u(\rho),\epsilon_0\}.$$
We also fix $\epsilon$ and $\eta$ two small positive constants lower than the injectivity radius of the manifold (that we suppose to be larger than $2$). We choose $\eta$ small enough to have $(2+\frac{b_0}{\epsilon_0})b_0\eta\leq\frac{\epsilon}{2}$ (as in~\cite{GR}, this property is only used in the proof of lemma~\ref{adaptedbis}). We underline that there exists $\varepsilon>0$ such that if
$$\forall\ (\rho,\rho')\in\mathcal{E}^{\theta}\times\mathcal{E}^{\theta},\ d(\rho,\rho')\leq\varepsilon\Rightarrow |U^u(\rho)-U^u(\rho')|\leq \epsilon_0\epsilon.$$
We make the extra assumption that the small parameter $\epsilon$ used for the continuity is smaller than $\epsilon_0$.

\subsubsection*{Discretization of the manifold}

As in the case of Anosov surfaces, our strategy to prove theorem~\ref{maintheo2} will be to introduce a discrete reparametrization of the geodesic flow. Regarding this goal, we cut the manifold $M$ and precisely,
we consider a partition $M=\bigsqcup_{i=1}^K O_i$ of diameter smaller than some positive $\delta$. Let $(\Omega_i)_{i=1}^K$ be a finite open cover of $M$ such that for all $1\leq i\leq K$, $O_i\subsetneq\Omega_i$. For $\gamma\in\{1,\cdots,K\}^2$, define an open subset of $T^*M$:
$$V_{\gamma}:=(T^*\Omega_{\gamma_0}\cap g^{-\eta}T^*\Omega_{\gamma_1})\cap\mathcal{E}^{\theta}.$$
We choose the partition $(O_i)_{i=1}^K$ and the open cover $(\Omega_i)_{i=1}^K$ of $M$ such that  $(V_{\gamma})_{\gamma\in\{1,\cdots,K\}^2}$ is a finite open cover of diameter smaller\footnote{In particular, the diameter of the partition $\delta$ depends on $\theta$ and $\epsilon$.} than $\varepsilon$ of $\mathcal{E}^{\theta}$. For $\gamma:=(\gamma_0,\gamma_1)$, we define $f(\gamma)$ and $f_{0}(\gamma)$ as in the case of an Anosov flow i.e.
$$f_0(\gamma):=\eta\inf\{U_0^u(\rho):\rho\in V_{\gamma}\}\ \text{and}\ f(\gamma):=\eta\inf\{U^u(\rho):\rho\in V_{\gamma}\}.$$
Compared with the Anosov case, we will have slightly different properties for the function $f(\gamma)$, i.e.
\begin{equation}\label{continuity2}\forall\rho\in V_{\gamma},\ \left|\int_{0}^{\eta}U^u_{\rho}(s)ds-f(\gamma)\right|\leq\eta\epsilon_0\epsilon.\end{equation}
We also underline that the function $f_0$ satisfies the following bounds, for $\gamma\in\{1,\cdots,K\}^2$,
$$\epsilon_0\eta\leq f_0(\gamma)\leq b_0\eta.$$
Finally, let $\alpha=(\alpha_0,\alpha_1,\cdots)$ be a sequence (finite or infinite) of elements of $\{1,\cdots,K\}$ whose length is larger than~$1$ and define:
\begin{equation}\label{ric}f_+(\alpha):=f_0\left(\alpha_0,\alpha_1\right)\leq\frac{\epsilon}{2}\ \text{and}\ f(\alpha):=f(\alpha_0,\alpha_1)\leq\frac{\epsilon}{2},\end{equation}
where the upper bounds follow from the previous hypothesis. In the following, we will also have to consider negative times. To do this, we define the analogous functions, for $\beta:=(\cdots,\beta_{-1},\beta_0)$ of finite (or infinite) length,
$$f_-(\beta):=f_0(\beta_{-1},\beta_0)\ \text{and}\ f(\beta):=f(\beta_{-1},\beta_0).$$
\begin{rema} We underline that the functions $f_+$ and $f_-$ are defined from $U^u_0$ while $f$ is defined from $U^u$. This distinction will be important in the following.
\end{rema}

\section{Proof of theorem~\ref{maintheo2}}\label{s:outline}

Let $(\psi_{\hbar_k})$ be a sequence of orthonormal eigenfunctions of the Laplacian corresponding to the eigenvalues $-1/\hbar_k^{-2}$ such that the corresponding sequence of distributions $\mu_k$ on $T^{*}M$ converges as $k$ tends to infinity to the semiclassical measure $\mu$. For simplicity of notations and to fit semiclassical analysis notations, we will denote $\hbar$ tends to $0$ the fact that $k$ tends to infinity and $\psi_{\hbar}$ and $\hbar^{-2}$ the corresponding eigenvector and eigenvalue. To prove the inequality of theorem~\ref{maintheo2}, we will again give a symbolic interpretation of a semiclassical measure and apply results on suspension flows to this measure~\cite{Ab}.\\
Let $\epsilon'>4\epsilon$ be a positive number, where $\epsilon$ was defined in section~\ref{classical-nonpos}. As in the Anosov setting, the link between the two quantities $\epsilon$ and $\epsilon'$ is only used to obtain theorem on product of pseudodifferential operators from sections $6$ and~$7$ in~\cite{GR} (here theorem~\ref{subadditivitypressure}). In the following of the note, the Ehrenfest time $n_E(\hbar)$ will be the quantity:
\begin{equation}\label{ehrenfest}n_E(\hbar):=[(1-\epsilon')|\log\hbar|].\end{equation}
We underline that it is an integer time and that, compared with usual definitions of the Ehrenfest time, there is no dependence on the Lyapunov exponent. We also consider a smaller non integer time:
\begin{equation}\label{ehrint}T_E(\hbar):=(1-\epsilon)n_E(\hbar).\end{equation}
We draw now a precise outline of the proof of theorem~\ref{maintheo2} and refer the reader to~\cite{GR} for the proof of several lemmas. The main differences with the Anosov case is that we have to indroduce a thermodynamical formalism to treat the problem.

\subsection{Quantum partitions of identity}

In order to find a lower bound on the metric entropy of the semiclassical measure $\mu$, we would like to apply the uncertainty principle for quantum pressure (see appendix~\ref{quant-pressure}) and see what informations it will give (when $\hbar$ tends to $0$) on the metric entropy of the semiclassical measure $\mu$. To do this, we define quantum partitions of identity corresponding to a given partition of the manifold. We recall the notations from~\cite{GR}.

\subsubsection{Partitions of identity}

In paragraph~\ref{classical-nonpos}, we considered a partition of small diameter $(O_i)_{i=1}^K$ of $M$. We also defined $(\Omega_i)_{i=1}^K$ a corresponding finite open cover of small diameter of $M$. By convolution of the characteristic functions $\mathbf{1}_{O_i}$, we obtain $\displaystyle\mathcal{P}=\left(P_i\right)_{i=1,..K}$ a smooth partition of unity on $M$ i.e. for all $x\in M$:
$$\sum_{i=1}^{K}P_i^2(x)=1.$$
We assume that for all $1\leq i\leq K$, $P_i$ is an element of $\mathcal{C}^{\infty}_c(\Omega_i)$. To this classical partition corresponds a quantum partition of identity of $L^2(M)$. In fact, if $P_i$ denotes the multiplication operator by $P_i(x)$ on $L^2(M)$, then one has:
\begin{equation}\label{time0bis}\sum_{i=1}^{K}P_i^{*}P_i=\text{Id}_{L^2(M)}.\end{equation}

\subsubsection{Refinement of the quantum partition under the Schrödinger flow}
\label{refinement0}
Like in the classical setting of entropy, we would like to make a refinement of the quantum partition. To do this refinement, we use the Schrödinger propagation operator $U^t=e^{\frac{\imath t\hbar\Delta}{2}}$. We define  $A(t):=U^{-t}AU^t$, where $A$ is an operator on $L^2(M)$. To fit as much as possible with the metric entropy, we define the following operators:
\begin{equation}\label{taubis}\tau_{\alpha}=P_{\alpha_{k}}(k\eta)\cdots P_{\alpha_{1}}(\eta)P_{\alpha_0}\end{equation}
and
\begin{equation}\label{pibis}\pi_{\beta}=P_{\beta_{-k}}(-k\eta)\cdots P_{\beta_{-2}}(-2\eta)P_{\beta_0}P_{\beta_{-1}}(-\eta),\end{equation}
where $\alpha=(\alpha_0,\cdots,\alpha_k)$ and $\beta=(\beta_{-k},\cdots,\beta_0)$ are finite sequences of symbols such that $\alpha_j\in [1,K]$ and $\beta_{-j}\in [1,K]$. We can remark that the definition of $\pi_{\beta}$ is the analogue for negative times of the definition of $\tau_{\alpha}$. The only difference is that we switch the two first terms $\beta_0$ and $\beta_{-1}$. The reason of this choice relies again in the application of the quantum uncertainty principle. One can see that for fixed $k$, using the Egorov property:
\begin{equation}\label{applegorovbis}\|P_{\alpha_{k}}(k\eta)\cdots P_{\alpha_{1}}(\eta)P_{\alpha_0}\psi_{\hbar}\|^2\rightarrow\mu(P_{\alpha_{k}}^2\circ g^{k\eta}\times\cdots P_{\alpha_{1}}^2\circ g^{\eta}\times P_{\alpha_0}^2)\ \text{as}\ \hbar\ \text{tends}\ \text{to}\ 0.\end{equation}
This last quantity is the one used to compute $h_{KS}(\mu,g^{\eta})$ (with the notable difference that the $P_j$ are here smooth functions instead of characteristic functions). As in~\cite{GR}, we will have to understand for which range of times $k\eta$, the Egorov property can be be applied. In particular, we will study for which range of times, the operator $\tau_{\alpha}$ is a pseudodifferential operator of symbol $P_{\alpha_{k}}\circ g^{k\eta}\times\cdots P_{\alpha_{1}}\circ g^{\eta}\times P_{\alpha_0}$ (see~(\ref{applegorovbis})). In~\cite{AN2} and~\cite{AKN}, they only considered $k\eta\leq|\log\hbar|/\lambda_{\max}$ where $\lambda_{\max}:=\lim_{t\rightarrow\pm\infty}\frac{1}{t}\log\sup_{\rho\in S^*M}|d_{\rho}g^t|$. This choice was not optimal and in the following, we try to define sequences $\alpha$ for which we can say that $\tau_{\alpha}$ is a pseudodifferential operator.

\subsubsection{Index family adapted to the variation of the unstable Jacobian}
\label{roofbis}
Let $\alpha=(\alpha_0,\alpha_1,\cdots)$ be a sequence (finite or infinite) of elements of $\{1,\cdots,K\}$ whose length is larger than~$1$. We define a natural shift on these sequences
$$\sigma_+((\alpha_0,\alpha_1,\cdots)):=(\alpha_1,\cdots).$$
For negative times and for $\beta:=(\cdots,\beta_{-1},\beta_0)$, we define the backward shift
$$\sigma_-((\cdots,\beta_{-1},\beta_0)):=(\cdots,\beta_{-1}).$$
In the following, we will mostly use the symbol $x$ for infinite sequences and reserve $\alpha$ and $\beta$ for finite ones. Then, using notations of section~\ref{nonpositive} and as described in section~$4$ of~\cite{GR}, index families depending on the value of the unstable Jacobian can be defined as follows:
\begin{equation}\label{Ihbis}I^{\eta}(\hbar):=I^{\eta}(T_E(\hbar))=\left\{\left(\alpha_0,\cdots,\alpha_{k}\right):k\geq3,\sum_{i=1}^{k-2}f_{+}\left(\sigma^i_+\alpha\right)\leq T_E(\hbar)<\sum_{i=1}^{k-1}f_{+}\left(\sigma^i_+\alpha\right)\right\},\end{equation}
\begin{equation}\label{Khbis}K^{\eta}(\hbar):=K^{\eta}(T_E(\hbar))=\left\{\left(\beta_{-k},\cdots,\beta_{0}\right): k\geq3,\sum_{i=1}^{k-2}f_{-}\left(\sigma_-^{i}\beta\right)\leq T_E(\hbar)<\sum_{i=1}^{k-1}f_{-}\left(\sigma_-^{i}\beta\right)\right\}.\end{equation}
We underline that $f_+,f_-\geq\epsilon_0\eta$ ensures that we consider finite sequences. These sets define the maximal sequences for which we can expect to have Egorov property for the corresponding $\tau_{\alpha}$. The sums used to define these sets were already used in~\cite{GR}. We can again think of the time $|\alpha|\eta$ as a stopping time for which property~(\ref{applegorovbis}) will hold (for a symbol $\tau_{\alpha}$ corresponding to $\alpha$).\\
A good way of thinking of these families of words is by introducing the sets
$$\Sigma_+:=\{1,\cdots,K\}^{\mathbb{N}}\ \text{and}\ \Sigma_-:=\{1,\cdots,K\}^{-\mathbb{N}}.$$
Once more, the sets $I^{\eta}(\hbar)$ (resp. $K^{\eta}(\hbar)$) lead to natural partitions of $\Sigma$ (resp. $\Sigma_-$). Families of operators can be associated to these families of index: $(\tau_{\alpha})_{\alpha\in I^{\eta}(\hbar)}$ and $(\pi_{\beta})_{\beta\in K^{\eta}(\hbar)}$. One can show that these partitions form quantum partitions of identity (lemma $5.1$ in~\cite{GR}):
$$\sum_{\alpha\in I^{\eta}(\hbar)}\tau_{\alpha}^*\tau_{\alpha}=\text{Id}_{L^{2}(M)}\ \text{and}\ \sum_{\beta\in K^{\eta}(\hbar)}\pi_{\beta}^*\pi_{\beta}=\text{Id}_{L^{2}(M)}.$$

\subsection{Symbolic interpretation of semiclassical measures}
\label{symbolicbis}
Now that we have defined these partitions of variable size, we want to show that they are adapted to compute the pressure of a certain measure with respect to some reparametrized flow associated to the geodesic flow. To do this, we proceed as in~\cite{GR} and provide a symbolic interpretation of the quantum partitions. We denote $\Sigma_+:=\{1,\cdots,K\}^{\mathbb{N}}$. We also denote $\mathcal{C}_i$ the subset of sequences $(x_n)_{n\in\mathbb{N}}$ such that $x_0=i$. Define also:
$$[\alpha_{0},\cdots,\alpha_{k}]:=\mathcal{C}_{\alpha_{0}}\cap\cdots\cap\sigma^{-k}_+\mathcal{C}_{\alpha_{k}},$$
where $\sigma_+$ is the shift $\sigma_+((x_n)_{n\in\mathbb{N}})=(x_{n+1})_{n\in\mathbb{N}}$ (it fits the notations of the previous section). The set $\Sigma_+$ is then endowed with the probability measure (not necessarily $\sigma$-invariant):
$$\mu_{\hbar}^{\Sigma_+}\left(\left[\alpha_{0},\cdots,\alpha_{k}\right]\right)=\mu_{\hbar}^{\Sigma_+}\left(\mathcal{C}_{\alpha_{0}}\cap\cdots\cap\sigma^{-k}_+\mathcal{C}_{\alpha_{k}}\right)=\|P_{\alpha_{k}}(k\eta)\cdots P_{\alpha_{0}}\psi_{\hbar}\|^2.$$
Using the property of partition of identity, it is clear that this definition assures the compatibility conditions to define a probability measure:
$$\sum_{\alpha_{k+1}}\mu_{\hbar}^{\Sigma_+}\left(\left[\alpha_{0},\cdots,\alpha_{k+1}\right]\right)=\mu_{\hbar}^{\Sigma_+}\left(\left[\alpha_{0},\cdots,\alpha_{k}\right]\right).$$
Then, we can define the suspension flow, in the sense of Abramov, associated to this probability measure. To do this, the suspension set is defined as:
\begin{equation}\label{suspsetbis}\overline{\Sigma}_+:=\{\left(x,s\right)\in\Sigma_+\times\mathbb{R}_+:0\leq s<f_+\left(x\right)\}.\end{equation}
Recall that the roof function $f_+$ is defined as $f_+(x):=f_+(x_0,x_1).$ We define a probability measure $\overline{\mu}_{\hbar}^{\overline{\Sigma}_+}$ on $\overline{\Sigma}_+$:
\begin{equation}\label{suspmeasbis}\overline{\mu}_{\hbar}^{\overline{\Sigma}_+}=\mu_{\hbar}^{\Sigma_+}\times \frac{dt}{\sum_{\alpha\in\{1,\cdots,K\}^2}f_+(\alpha)\|P_{\alpha}\psi_{\hbar}\|^2}=\mu_{\hbar}^{\Sigma_+}\times \frac{dt}{\sum_{\alpha\in\{1,\cdots,K\}^2}f_+(\alpha)\mu_{\hbar}^{\Sigma_+}\left(\left[\alpha\right]\right)}.\end{equation}
The suspension semi-flow associated to $\sigma_+$ is for time $s$:
\begin{equation}\label{suspflowbis}\overline{\sigma}^s_+\left(x,t\right):=\left(\sigma^{n-1}_+ (x),s+t-\sum_{j=0}^{n-2}f_+\left(\sigma^j_+x\right)\right),\end{equation}
where $n$ is the only integer such that $\displaystyle\sum_{j=0}^{n-2}f_+\left(\sigma^j_+x\right)\leq s+t<\sum_{j=0}^{n-1}f_+\left(\sigma^j_+x\right)$.
\begin{rema} We underline that we used the fact that $f_+>0$ to define the suspension flow. If we had considered $f$, we would not have been able to construct the suspension flow as $f$ could be equal to $0$.
\end{rema}
A notable difference with the Anosov setting is that we will not consider time $1$ of the suspension of the flow. Instead of it, we fix a large integer $N_0$ (such that\footnote{To summarize the relations between the different parameters, we have $\frac{\epsilon}{4}<\epsilon'\ll \frac{1}{N_0}\ll\epsilon_0$. Moreover $\eta$ depends on $\epsilon$ and $\epsilon_0$ and tends to $0$ when $\epsilon$ tends to $0$ and $\epsilon_0$ is fixed.} $\epsilon'\ll 1/N_0\ll \epsilon_0$) and consider time $1/N_0$ of the flow and its iterates.
\begin{rema} It can be underlined that the same procedure holds for the partition $(\pi_{\beta})$. The only differences are that we have to consider $\Sigma_{-}:=\{1,\cdots,K\}^{-\mathbb{N}}$, $\sigma_-((x_n)_{n\leq 0})=(x_{n-1})_{n\leq 0}$ and that the corresponding measure is, for $k\geq 1$:
$$\mu_{\hbar}^{\Sigma_{-}}\left(\left[\beta_{-k},\cdots,\beta_{0}\right]\right)=\mu_{\hbar}^{\Sigma_{-}}\left(\sigma_-^{-k}\mathcal{C}_{\beta_{-k}}\cap\cdots\cap\mathcal{C}_{\beta_{0}}\right)=\|P_{\beta_{-k}}(-k\eta)\cdots P_{\beta_{0}}P_{\beta_{-1}}(-\eta)\psi_{\hbar}\|^2.$$
For $k=0$, one should take the only possibility to assure the compatibility condition:
$$\mu_{\hbar}^{\Sigma_{-}}\left(\left[\beta_{0}\right]\right)=\sum_{j=1}^K\mu_{\hbar}^{\Sigma_{-}}\left(\left[\beta_{-1},\beta_{0}\right]\right).$$
The definition is quite different from the positive case but in the semiclassical limit, it will not change anything as $P_{\beta_0}$ and $P_{\beta_-1}(-\eta)$ commute. Finally, the "past" suspension set can be defined as
$$\overline{\Sigma}_-:=\{(x,s)\in\Sigma_-\times\mathbb{R}_+:0\leq s<f_-(x)\}.$$
\end{rema}
Now let $\alpha$ be an element of $I^{\eta}(\hbar)$. Define:
\begin{equation}\label{symbpart}\tilde{\mathcal{C}}_{\alpha}^+:=\mathcal{C}_{\alpha_0}\cap\cdots\cap\sigma^{-k}_+\mathcal{C}_{\alpha_{k}}.\end{equation}
This new family of subsets forms a partition of $\Sigma_+$. Then, a partition $\overline{\mathcal{C}}_{\hbar}^+$ of $\overline{\Sigma}_+$ can be defined starting from the partition $\tilde{\mathcal{C}}$ and $[0,f_+(\alpha)[$. An atom of this suspension partition is an element of the form $\overline{\mathcal{C}}_{\alpha}^+=\tilde{\mathcal{C}}_{\alpha}^+\times [0,f_+(\alpha)[$. For $\overline{\Sigma}^-$ (the suspension set corresponding to $\Sigma_-$), we define an analogous partition $\overline{\mathcal{C}}^{-}_{\hbar}=([\beta]\times[0,f_-(\beta)[)_{\beta\in K^{\eta}(\hbar)}$. As in the case of the Anosov geodesic flows, we now have to apply the uncertainty principle to these partitions of variable size. The main difference with~\cite{GR} is that we will apply it for quantum pressures (see section~\ref{quant-pressure}). We introduce the weights
$$W_{\alpha}^+:=\exp\left(\frac{1}{2}\sum_{j=1}^{k-1}f(\sigma^j_+\alpha)\right)\ \text{and}\ W_{\beta}^-:=\exp\left(\frac{1}{2}\sum_{j=1}^{k-1}f(\sigma^j_-\beta)\right).$$
We underline that the weights depends on $f$ and not $f_+$ or $f_-$. It cames from the fact that $f$ is the function that appears in theorem~\ref{normestbis}. We introduce the associated quantum pressure\footnote{We refer the reader to appendix~\ref{KSentropy} for the definition of $H$.}:
\begin{equation}\label{pressure}p\left(\overline{\mu}_{\hbar}^{\overline{\Sigma}_+},\overline{\mathcal{C}}_{\hbar}^+\right):=H\left(\overline{\mu}_{\hbar}^{\overline{\Sigma}_+},\overline{\mathcal{C}}_{\hbar}^+\right)-2\sum_{\alpha\in I^{\eta}(\hbar)}\overline{\mu}_{\hbar}^{\overline{\Sigma}_+}\left(\overline{\mathcal{C}}_{\alpha}^+\right)\log W_{\alpha}^+\end{equation}
and
\begin{equation}\label{pressureneg}p\left(\overline{\mu}_{\hbar}^{\overline{\Sigma}_-},\overline{\mathcal{C}}_{\hbar}^-\right):=H\left(\overline{\mu}_{\hbar}^{\overline{\Sigma}_-},\overline{\mathcal{C}}_{\hbar}^-\right)-2\sum_{\beta\in K^{\eta}(\hbar)}\overline{\mu}_{\hbar}^{\overline{\Sigma}_-}\left(\overline{\mathcal{C}}_{\beta}^-\right)\log W_{\beta}^-.\end{equation}
We follow then the procedure of section $5$ in~\cite{GR} to apply the entropic uncertainty principle (i.e. apply it $K^2$ times and not $1$ time as in~\cite{AKN}) and we use the main estimate on the norms of the quantum partitions (see theorem~\ref{normestbis}) to derive that
\begin{equation}\label{mainest-nonpos}p\left(\overline{\mu}_{\hbar}^{\overline{\Sigma}_+},\overline{\mathcal{C}}_{\hbar}^+\right)+p\left(\overline{\mu}_{\hbar}^{\overline{\Sigma}_-},\overline{\mathcal{C}}_{\hbar}^-\right)\geq -\log C-(1+\epsilon'+4\epsilon)n_E(\hbar),\end{equation}
where $C$ is a constant that does not depend on $\hbar$.
\begin{rema} This last inequality is a crucial step to prove theorem~\ref{maintheo2}. We will recall how one can get such a lower bound in section~\ref{s:appl-uncert}. This inequality corresponds to proposition $5.3$ in~\cite{GR}. The strategy of the proof is exactly the same except that we have to deal with quantum pressures and not quantum entropies (see section~\ref{s:appl-uncert}). However, we can follow the same lines as in section $5.3.2$ in~\cite{GR} (i.e. apply $K^2$ times the uncertainty principle) and obtain a lower bound that depends on the bound from theorem~\ref{normestbis}. At this point, there is a difference because theorem~\ref{normestbis} was proved in~\cite{AN2} for Anosov manifolds. In section~\ref{s:main-est}, we will show that the proof of theorem~\ref{normestbis} from~\cite{AN2} can be adapted in the setting of nonpositively curved surfaces.
\end{rema}
The problem of expression~(\ref{mainest-nonpos}) is that it is not exactly the pressure of a refined partition. As in~\cite{GR}, one can prove the following lemma:
\begin{lemm}\label{adaptedbis} Let $N_0$ be a positive integer defined as previously. There exists an explicit partition $\overline{\mathcal{C}}^+_{N_0}$ of $\overline{\Sigma}_+$, independent of $\hbar$ such that $\vee_{i=0}^{n_E(\hbar)N_0-1}\overline{\sigma}^{-\frac{i}{N_0}}_+\overline{\mathcal{C}}_+$ is a refinement of the partition $\overline{\mathcal{C}}_{\hbar}^+$. Moreover, let $n$ be a fixed positive integer. Then, an atom of the refined partition $\displaystyle\vee_{i=0}^{n-1}\overline{\sigma}^{-\frac{i}{N_0}}_+\overline{\mathcal{C}}_+$ is of the form $[\alpha]\times B(\alpha)$, where $\alpha=(\alpha_0,\cdots,\alpha_k)$ is a $k+1$-uple such that $(\alpha_0,\cdots,\alpha_k)$ verifies $\frac{n}{N_0}(1-\epsilon)\leq\displaystyle\sum_{j=0}^{k-1}f_+\left(\sigma^j_+\alpha\right)\leq \frac{n}{N_0}(1+\epsilon)$ and $B(\alpha)$ is a subinterval of $[0,f_+(\alpha)[$.
\end{lemm}
This lemma is the exact analogue of lemma $4.1$ in~\cite{GR} and its proof is the same: the only difference is that we consider times $1/N_0$ instead of time $1$. In particular, in the proof, the partition $\overline{\mathcal{C}}_{N_0}^+$ is constructed from\footnote{We recall that $I^{\eta}(t)$ was defined as the set of words $\left\{\alpha=\left(\alpha_0,\cdots,\alpha_{k}\right):k\geq3,\sum_{i=1}^{k-2}f_+\left(\sigma^i_+\alpha\right)\leq t<\sum_{i=1}^{k-1}f_+\left(\sigma^i_+\alpha\right)\right\}$.} $I^{\eta}(1/N_0)$ and not from $I^{\eta}(1)$. As in the Anosov case, we would like to use this lemma to rewrite the quantum pressure in terms of the pressure of a refined partition. To do this, we use basic properties of the classical entropy (see appendix~\ref{KSentropy}) to find that:
$$H\left(\overline{\mu}_{\hbar}^{\overline{\Sigma}_+},\overline{\mathcal{C}}_{\hbar}^+\right)\leq H_{N_0n_E(\hbar)}\left(\overline{\mu}_{\hbar}^{\overline{\Sigma}_+},\overline{\sigma}^{\frac{1}{N_0}}_+,\overline{\mathcal{C}}_{N_0}^+\right).$$
Consider now an atom $A$ of the partition $\vee_{j=0}^{n_{E}(\hbar)N_0-1}\overline{\sigma}_+^{-\frac{j}{N_0}}\overline{\mathcal{C}}_{N_0}^+$. To this atom, it corresponds an unique family $(\gamma_0,\cdots,\gamma_{n_E(\hbar)N_0-1})$ in $I^{\eta}(1/N_0)^{N_0n_E(\hbar)}$ and we define the corresponding weight as
$$W_A^+:=\prod_{j=0}^{N_0n_E(\hbar)-1}W_{\gamma_j}^+.$$
With these notations, we introduce the refined pressure at times $n$:
$$p_{n}\left(\overline{\mu}_{\hbar}^{\overline{\Sigma}_+},\overline{\sigma}^{\frac{1}{N_0}}_+,\overline{\mathcal{C}}^+_{N_0}\right):=H_{n}\left(\overline{\mu}_{\hbar}^{\overline{\Sigma}_+},\overline{\sigma}^{\frac{1}{N_0}}_+,\overline{\mathcal{C}}^+_{N_0}\right)-2\sum_{A\in\vee_{j=0}^{n-1}\overline{\sigma}^{-\frac{j}{N_0}}\overline{\mathcal{C}}_{N_0}^+}\overline{\mu}_{\hbar}^{\overline{\Sigma}_+}(A)\log W_A^+.$$
One can then write the following inequality
$$-2\sum_{\alpha\in I^{\eta}(\hbar)}\overline{\mu}_{\hbar}^{\overline{\Sigma}_+}\left(\overline{\mathcal{C}}_{\alpha}^+\right)\log W_{\alpha}^+\leq-2\sum_{A\in\vee_{j=0}^{N_0n_{E}(\hbar)-1}\overline{\sigma}^{-\frac{j}{N_0}}\overline{\mathcal{C}}_{N_0}^+}\overline{\mu}_{\hbar}^{\overline{\Sigma}_+}(A)\log W_A^++2\frac{b_0}{\epsilon_0}b_0\eta N_0n_E(\hbar).$$
The correction term in the last expression comes from the fact that, for each atom $A$ in the partition $\vee_{j=0}^{N_0n_{E}(\hbar)-1}\overline{\sigma}^{-\frac{j}{N_0}}\overline{\mathcal{C}}$, one has an unique $\alpha'$ in $I^{\eta}(\hbar)$ and the corresponding $W_{\alpha'}^+$ is not exactly equal to $W_A^+$. Finally, the previous inequalities can be summarized as follows:
\begin{equation}\label{crucial-est}-4\frac{b_0N_0n_E(\hbar)}{\epsilon_0}\epsilon-\log C-(1+\epsilon'+4\epsilon)n_E(\hbar)\leq p_{n_E(\hbar)N_0}\left(\overline{\mu}_{\hbar}^{\overline{\Sigma}_+},\overline{\sigma}^{\frac{1}{N_0}}_+,\overline{\mathcal{C}}^+_{N_0}\right)+p_{n_E(\hbar)N_0}\left(\overline{\mu}_{\hbar}^{\overline{\Sigma}_-},\overline{\sigma}^{\frac{1}{N_0}}_-,\overline{\mathcal{C}}^-_{N_0}\right).\end{equation}
This estimate is crucial in our proof as we have derived from a quantum relation a lower bound on the classical pressure of a dynamical system associated to the geodesic flow.

\subsection{Subadditivity of the quantum pressure}

As in~\cite{GR}, we would like to let $\hbar$ tends to $0$ in inequality~(\ref{crucial-est}). The main difficulty to do this is that everything depends on $\hbar$. So, once more, we have to prove a subadditivity property for the quantum pressure:
\begin{theo} \label{subadditivitypressure}Let $\overline{\mathcal{C}}_{N_0}^+$ be the partition of lemma~\ref{adaptedbis}. There exists a function $R(n_0,\hbar)$ on $\mathbb{N}\times(0,1]$ and $R(N_0)$ independent of $n_0$ such that
$$\forall n_0\in\mathbb{N},\ \ \ \ \limsup_{\hbar\rightarrow 0}|R(n_0,\hbar)|=R(N_0).$$
Moreover, for any $\hbar\in(0,1]$ and any $n_0,m\in\mathbb{N}$ such that $n_0+m\leq N_0n_E(\hbar)$, one has:
$$p_{n_0+m}\left(\overline{\mu}_{\hbar}^{\overline{\Sigma}_+},\overline{\sigma}^{\frac{1}{N_0}},\overline{\mathcal{C}}_{N_0}^+\right)\leq p_{n_0}\left(\overline{\mu}_{\hbar}^{\overline{\Sigma}_+},\overline{\sigma}^{\frac{1}{N_0}},\overline{\mathcal{C}}_{N_0}^+\right)+p_{m}\left(\overline{\mu}_{\hbar}^{\overline{\Sigma}_+},\overline{\sigma}^{\frac{1}{N_0}},\overline{\mathcal{C}}_{N_0}^+\right)+R(n_0,\hbar).$$
\end{theo}
\emph{Proof.} To prove this subadditivity property, we will prove subadditivity of the quantum entropy and subadditivity of the pressure term. As in section $6$ from~\cite{GR}, we write for the entropy part that:
$$H_{n_0+m}\left(\overline{\mu}_{\hbar}^{\overline{\Sigma}_+},\overline{\sigma}^{\frac{1}{N_0}},\overline{\mathcal{C}}_{N_0}^+\right)\leq H\left(\overline{\sigma}_+^{\frac{m}{N_0}}\sharp\overline{\mu}_{\hbar}^{\overline{\Sigma}_+},\vee_{j=0}^{n_0-1}\overline{\sigma}^{-\frac{j}{N_0}}\overline{\mathcal{C}}_{N_0}^+\right)+H_{m}\left(\overline{\mu}_{\hbar}^{\overline{\Sigma}_+},\overline{\sigma}^{\frac{1}{N_0}},\overline{\mathcal{C}}_{N_0}^+\right).$$
So, as in~\cite{GR}, we have to show that the measure of the atoms of the partition are almost invariant under $\overline{\sigma}_+^{\frac{1}{N_0}}$ for the range of times we have considered (proposition~$6.1$ in~\cite{GR}). Consider now the pressure term in the quantum pressure. Using the multiplicative structure of the $W_A^+$, one has
$$\sum_{A\in\vee_{j=0}^{n_0+m-1}\overline{\sigma}^{-\frac{j}{N_0}}\overline{\mathcal{C}}_{N_0}^+}\overline{\mu}_{\hbar}^{\overline{\Sigma}_+}(A)\log W_A^+=\sum_{A\in\vee_{j=0}^{m-1}\overline{\sigma}^{-\frac{j}{N_0}}\overline{\mathcal{C}}_{N_0}^+}\overline{\mu}_{\hbar}^{\overline{\Sigma}_+}(A)\log W_A^+\ \ \ \ \ \ \ \ \ \ \ \ \ \ \ \ \ \ \ \ $$
$$\ \ \ \ \ \ \ \ \ \ \ \ \ \ \ \ \ \ \ \ \ \ \ \ \ +\sum_{A\in\vee_{j=0}^{n_0-1}\overline{\sigma}^{-\frac{j}{N_0}}\overline{\mathcal{C}}_{N_0}^+}\overline{\sigma}_+^{\frac{m}{N_0}}\sharp\overline{\mu}_{\hbar}^{\overline{\Sigma}_+}(A)\log W_A^++\sum_{A\in\overline{\mathcal{C}}_{N_0}^+}\overline{\sigma}_+^{\frac{m}{N_0}}\sharp\overline{\mu}_{\hbar}^{\overline{\Sigma}_+}(A)\log W_A^+.$$
So, once more, the additivity property of the pressure term derives from the almost invariance of the measure for the range of times we consider\footnote{We underline that $R(N_0)$ will be equal to $\sup_{A\in\overline{\mathcal{C}}_{N_0}^+}\log W_A^+$ which only depends on $N_0$.}. Precisely, according to the last two inequalities, we only need to verify that proposition $6.1$ in~\cite{GR} remains true for the partition $\overline{\mathcal{C}}_{N_0}^+$ in the setting of surfaces of nonpositive curvature. We will not reproduce the proof here: it is the same one. We recall that this proposition relied on a theorem for products of pseudodifferential operators (theorem $7.1$ in~\cite{GR}) and we need to verify that the proof we gave still works in the case of surfaces of nonpositive curvature. The key point of the proof of this theorem is that in the allowed range of times $|d_{\rho}g^t|$ is bounded by some $\hbar^{-\nu}$ (with $\nu<1/2$) (see section~$7.2$ in~\cite{GR}). We know that to each $\rho$ we can associate a word $\alpha$ of length $k$. The range of times we will consider will be $0\leq t\leq k\eta$. To prove previous property in the case of surfaces of nonpositive curvature, we use the splitting of $T_{\rho}S^*M$ given by $\mathbb{R}X_H(\rho)\oplus E^u(\rho)\oplus V_{\rho}$. These three subspaces are uniformly transverse so we only have to give an estimate of $\|d_{\rho}g^t_{E\rightarrow T^*_{g^t\rho}M}\|$ when $E$ is one of them. In the case where $E=\mathbb{R}X_H(\rho)$, it is bounded by $1$ and in the case where $E=E^u(\rho)$, it is bounded by $\sqrt{1+K_0}e^{\int_0^{t}U^u_{\rho}(s)ds}$. In the last case, lemma~\ref{unstable} tells us that the spaces $d_{\rho}g^tV_{\rho}$ and $E^u(g^t\rho)$ become uniformly close (in direction) to each other. Then, we consider $e_0$ a unit vector in $V_{\rho}$ and for $0\leq p\leq k-1$, we define the $e_{p\eta}$ as the unit vector $\frac{d_{\rho}g^{p\eta}e_0}{\|d_{\rho}g^{p\eta}e_0\|}$. We can write:
$$\|d_{\rho}g^{k\eta}e_0\|=|\langle d_{\rho}g^{k\eta}e_0,e_{k\eta}\rangle|=|\langle d_{g^{(k-1)\eta}\rho}g^{\eta}e_{(k-1)\eta},e_{k\eta}\rangle\cdots \langle d_{\rho}g^{\eta}e_{0},e_{\eta}\rangle|.$$
We also define the corresponding sequence $e_{p\eta}^u:=\frac{d_{\rho}g^{p\eta}e_0^u}{\|d_{\rho}g^{p\eta}e_0^u\|}$ of unit unstable vectors, where $e_0^u:=\frac{(J^u_{\rho}(0),J^{u'}_{\rho}(0))}{\|(J^u_{\rho}(0),J^{u'}_{\rho}(0))\|}$. From lemma~\ref{unstable}, we know that $e_{p\eta}$ becomes uniformly close (in $\rho$) to $e_{p\eta}^u$. So, up to an error term of order $Ce^{k\eta\delta}$ (with $C$ uniform in $\rho$ and $\delta$ arbitrarly small), we have:
$$\|d_{\rho}g^{k\eta}e_0\|\leq Ce^{k\eta\epsilon_0\epsilon}|\langle d_{g^{(k-1)\eta}\rho}g^{\eta}e_{(k-1)\eta}^u,e_{k\eta}^u\rangle\cdots \langle d_{\rho}g^{\eta}e_{0}^u,e_{\eta}^u\rangle|=Ce^{k\eta\delta}\|d_{\rho}g^{k\eta}_{|E^u(\rho)}\|.$$
Finally, taking $\delta=\epsilon_0\epsilon$, we have that $\|d_{\rho}g^{k\eta}\|$ is bounded by $Ce^{k\eta\epsilon_0\epsilon}e^{\int_0^{k\eta}U^u_{\rho}(s)ds}$ (with $C$ uniform in $\rho$). For the allowed words, $e^{k\eta\epsilon_0\epsilon}$ is of order $\hbar^{-\epsilon}$ (as $k\eta\epsilon_0\leq1/2n_E(\hbar)$). To conclude, we can estimate:
$$\left|\int_0^{k\eta}U^u_{\rho}(s)ds-\sum_{j=0}^{k-1}f(\sigma^j\alpha)\right|\leq\sum_{j=0}^{k-1}\left|\int_{j\eta}^{(j+1)\eta}U^u_{\rho}(s)ds-f(\sigma^j\alpha)\right|.$$
To bound this sum, we can use the continuity of $U^u$ (see inequality~(\ref{continuity2})) to show that this quantity is bounded by $\epsilon|\log\hbar|$. By definition of the allowed words $\alpha$, we know that $\sum_{j=0}^{k-1}f(\sigma^j\alpha)\leq 1/2n_E(\hbar)$. This allows to conclude that $|d_{\rho}g^t|$ is bounded by some $C\hbar^{-\nu}$ (with $C$ independent of $\rho$ and $\nu<1/2$).$\square$

\begin{rema} We underline that here we need to use the specific properties of surfaces of nonpositive curvature to prove this theorem. It is not really surprising that theorem $7.1$ from~\cite{GR} can be extended in our setting as the situation can only be less `chaotic'. We also mention that we have to use the continuity of $U^u(\rho)$ which is for instance false for surfaces without conjugate points~\cite{BBB}.
\end{rema}

\subsection{The conclusion}

\subsubsection{Applying the Abramov theorem}
Thanks to the subadditivity property of the quantum pressure, we can proceed as in~\cite{GR} and write, for
a fixed $n_0$, the euclidean division $N_0n_E(\hbar)=qn_0+r$. Using the same method, we find, after applying the subadditivity property and letting $\hbar$ tends to $0$,
$$-4\frac{b_0}{\epsilon_0}\epsilon-\frac{R(N_0)}{n_0}-\frac{1}{N_0}(1+\epsilon'+4\epsilon)\leq \frac{1}{n_0}\left(p_{n_0}\left(\overline{\mu}^{\overline{\Sigma}_+},\overline{\sigma}^{\frac{1}{N_0}}_+,\overline{\mathcal{C}}^+_{N_0}\right)+p_{n_0}\left(\overline{\mu}^{\overline{\Sigma}_-},\overline{\sigma}^{\frac{1}{N_0}}_-,\overline{\mathcal{C}}^-_{N_0}\right)\right).$$
As in~\cite{GR}, we can replace the smooth partitions by true partitions of the manifold in the previous inequality. We would like now to transform the previous inequality on the metric pressure into an inequality on the Kolmogorov-Sinai entropy. To do this, we write the multiplicative property of $W_A$ to write:
$$\sum_{A\in\vee_{j=0}^{n_0-1}\overline{\sigma}^{-\frac{j}{N_0}_+}\overline{\mathcal{C}}^+_{N_0}}\overline{\mu}^{\overline{\Sigma}^+}(A)\log W_A^+=\sum_{A_0,\cdots,A_{n_0-1}\in\overline{\mathcal{C}}^+_{N_0}}\overline{\mu}^{\overline{\Sigma}^+}(A_0\cap\cdots\cap\overline{\sigma}^{-\frac{n_0-1}{N_0}}A_{n_0-1})\sum_{j=0}^{n_0-1}\log W_{A_j}^+.$$
After simplification and using the fact $\overline{\mathcal{C}}^+_{N_0}$ is a partition of $\overline{\Sigma}^+$, we find that this last inequality can be rewritten as follows
$$\sum_{A\in\vee_{j=0}^{n_0-1}\overline{\sigma}^{-\frac{j}{N_0}_+}\overline{\mathcal{C}}^+_{N_0}}\overline{\mu}^{\overline{\Sigma}^+}(A)\log W_A^+=n_0\sum_{A\in\overline{\mathcal{C}}^+_{N_0}}\overline{\mu}^{\overline{\Sigma}^+}(A)\log W_A^+$$
The same property holds for the backward side. After letting $n_0$ tends to infinity, we find that:
$$-4\frac{b_0}{\epsilon_0}\epsilon-\frac{1}{N_0}(1+\epsilon'+4\epsilon)+2\left(\sum_{A\in\overline{\mathcal{C}}^+_{N_0}}\overline{\mu}^{\overline{\Sigma}^+}(A)\log W_A^++\sum_{A\in\overline{\mathcal{C}}^-_{N_0}}\overline{\mu}^{\overline{\Sigma}^-}(A)\log W_A^-\right)$$
$$\ \ \ \ \ \ \ \ \ \ \ \ \ \ \ \ \ \ \ \ \ \ \ \ \ \ \ \ \ \ \ \ \ \ \ \ \leq \frac{1}{N_0}\left(h_{KS}\left(\overline{\mu}^{\overline{\Sigma}_+},\overline{\sigma}_+\right)+h_{KS}\left(\overline{\mu}^{\overline{\Sigma}_-},\overline{\sigma}_-\right)\right).$$
We now underline that, by construction (see the proof of lemma~$4.1$ in~\cite{GR}) and by invariance of the measure $\mu^{\Sigma_+}$, one has:
$$\sum_{A\in\overline{\mathcal{C}}^+_{N_0}}\overline{\mu}^{\overline{\Sigma}^+}(A)\log W_A^++\sum_{A\in\overline{\mathcal{C}}^-_{N_0}}\overline{\mu}^{\overline{\Sigma}^-}(A)\log W_A^-=\frac{2}{\sum_{\gamma'\in\{1,\cdots,K\}^2}f_0(\gamma')\mu^{\Sigma}([\gamma'])}\sum_{\gamma\in I^{\eta}(1/N_0)}f_0(\gamma)\mu^{\Sigma}([\gamma])\log W_{\gamma}.$$
We use this last property and combine it with the Abramov theorem~\cite{Ab}. We find then
$$\sum_{\gamma'\in\{1,\cdots,K\}^2}f_0(\gamma')\mu^{\Sigma}([\gamma'])\left(-2\frac{b_0N_0}{\epsilon_0}\epsilon-\frac{1}{2}(1+\epsilon'+4\epsilon)\right)+2N_0\sum_{\gamma\in I^{\eta}(1/N_0)}f_0(\gamma)\mu^{\Sigma}([\gamma])\log W_{\gamma}\leq \eta h_{KS}(\mu,g).$$

\subsubsection{The different small parameters tend to $0$} We have obtained a lower bound on the Kolmogorv-Sinai entropy of the measure $\mu$. This lower bound depend on several small parameters that are linked to each other in the following way:
$$\epsilon<4\epsilon'\ll\frac{1}{N_0}\ll\epsilon_0.$$
Moreover the small parameter $\eta$ depends on $\epsilon$ and $\epsilon_0$. For a fixed $\epsilon_0$, it tends to $0$ when $\epsilon$ tends to $0$. We have now to be careful to transform our lower bound on the entropy of $\mu$ into the expected lower bound. To do this, we use the notations of section~\ref{nonpositive} and introduce, for $\rho\in S^*M$, the application
$$F_0(\rho):=\sum_{\gamma\in I^{\eta}(1/N_0)}f_0(\gamma)\log W_{\gamma} \mathbf{1}_{O_{\gamma_0}}(\rho)\cdots\mathbf{1}_{O_{\gamma_k}}\circ g^{k\eta}(\rho).$$
We underline that for each $\rho$ in $S^*M$, there exists an unique $\gamma$ in $I^{\eta}(1/N_0)$ such that $\mathbf{1}_{O_{\gamma_0}}(\rho)\cdots\mathbf{1}_{O_{\gamma_k}}\circ g^{k\eta}(\rho)$ is non zero (it is then equal to $1$). With this new function, the lower bound on the Kolmogorov-Sinai entropy can be rewritten as follows:
$$\sum_{\gamma'\in\{1,\cdots,K\}^2}f_0(\gamma')\mu^{\Sigma}([\gamma'])\left(-2\frac{b_0N_0}{\epsilon_0}\epsilon-\frac{1}{2}(1+\epsilon'+4\epsilon)\right)+2N_0\int_{S^*M}F_0(\rho)d\mu(\rho)\leq \eta h_{KS}(\mu,g).$$
We define then
$$X_0:=\left\{\rho\in S^*M:\forall  0\leq t\leq\frac{1}{N_0\epsilon_0}, U^u(g^{t}\rho)>2\epsilon_0\right\}.$$
We can verify that $F_0(\rho)\geq (1/N_0)\sum_{\gamma_0,\gamma_1}f_0(\gamma)\mathbf{1}_{X_0}(\rho)\mathbf{1}_{O_{\gamma_0}}(\rho)\mathbf{1}_{O_{\gamma_1}}\circ g^{\eta}(\rho)$ for all $\rho$ in $\mathcal{E}^{\theta}$. In fact, one has, for $\rho\in X_0$ (otherwise the inequality is trivial), $\log W_{\gamma}=\frac{1}{2}\sum_{j=1}^{k-1}f(\sigma^j\gamma)$, where $\rho$ belongs to $O_{\gamma_0}\cap\cdots g^{-k\eta}O_{\gamma_k}$ and $\gamma$ satisfies
$$\sum_{j=1}^{k-2}f_0(\sigma^j\gamma)\leq\frac{1}{N_0}<\sum_{j=1}^{k-1}f_0(\sigma^j\gamma).$$
In particular, one has $(k-2)\eta\epsilon_0\leq 1/N_0$. Using the relation of continuity~(\ref{continuity2}) and the fact that $U^u_0(g^t\rho)=U^u(g^t\rho)$ on $X_0$, one find that, for $\rho\in X_0\cap O_{\gamma_0}\cap\cdots g^{-k\eta}O_{\gamma_k}$,
$$\log W_{\gamma}\geq-\frac{2\epsilon_0}{N_0}+\frac{1}{2}\sum_{j=1}^{k-1}f_0(\sigma^j\gamma)\geq\left(\frac{1}{2}-2\epsilon_0\right)\frac{1}{N_0}.$$
We use this function $\mathbf{1}_{X_0}(\rho)$ in our lower bound on the entropy of $\mu$. We let the diameter of the partition tends to $0$ and we divide by $\eta$. This gives us
$$\left(-2\frac{b_0N_0}{\epsilon_0}\epsilon-\frac{1}{2}(1+\epsilon'+4\epsilon)\right)\int_{S^*M}U^u_0(\rho)d\mu(\rho)+(1-4\epsilon_0)\int_{S^{*}M}U^u_0(\rho)\mathbf{1}_{X_0}(\rho)d\mu(\rho)\leq h_{KS}(\mu,g).$$
Finally, we let $\epsilon$ and $\epsilon'$ tend to $0$ (in this order). We find the following bound on the entropy of $\mu$:
$$-\frac{1}{2}\int_{S^*M}U^u_0(\rho)d\mu(\rho)+(1-4\epsilon_0)\int_{S^{*}M}U^u_0(\rho)\mathbf{1}_{X_0}(\rho)d\mu(\rho)\leq h_{KS}(\mu,g).$$
We let now $N_0$ tend to infinity and then $\epsilon_0$ tend to $0$ (in this order). We find finally the expected lower bound:
$$\frac{1}{2}\int_{S^*M}U^u(\rho)d\mu(\rho)\leq h_{KS}(\mu,g).\square$$

\section{Proof of the main estimate from~\cite{AN2}}\label{s:main-est}

In the previous section, we have been able to apply the method we used for Anosov surfaces in order to prove theorem~\ref{maintheo2}. As in~\cite{GR}, the strategy relied on a careful adaptation of an uncertainty principle. In particular, to derive inequality~(\ref{mainest-nonpos}), we had to use the following equivalent of theorem~$3.1$ from~\cite{AKN}:
\begin{theo}\label{normestbis} Let $M$ be a surface of nonpositive sectional curvature and $\epsilon$, $\epsilon_0$ and $\eta$ be small positive parameters as in section~\ref{classical-nonpos}. For every $\mathcal{K}>0$ ($\mathcal{K}\leq C_{\delta_0}$), there exists $\hbar_{\mathcal{K}}$ and $C_{\mathcal{K}}(\epsilon,\eta,\epsilon_0)$ such that uniformly for all $\hbar\leq\hbar_{\mathcal{K}}$, for all $k\leq\mathcal{K}|\log\hbar|$, for all $\alpha=(\alpha_0,\cdots,\alpha_k)$,
\begin{equation}\label{ineqnormbis}\|P_{\alpha_{k}}U^{\eta}P_{\alpha_{k-1}}\cdots U^{\eta}P_{\alpha_0}\Op_{\hbar}(\chi^{(k)})\|_{L^2(M)} \leq C_{\mathcal{K}}(\epsilon,\eta,\epsilon_0)\hbar^{-\frac{1}{2}-c\delta_0}e^{2k\eta\epsilon_0\epsilon}\exp\left(-\frac{1}{2}\sum_{j=0}^{k-1}f(\sigma^j_+\alpha)\right),\end{equation}
where $c$ depends only on the riemannian manifold $M$.
\end{theo}
\begin{rema} We underline two facts about this theorem. The first one is that $\Op_{\hbar}(\chi^{(k)})$ is a cutoff operator that was already defined in~\cite{GR} (section $5.3$) and in the appendix of~\cite{AN2}. We describe briefly its construction in section~\ref{cut}. The second one is that it is function $f$ and not $f_+$ that appears in the upper bound.
\end{rema}
This theorem is the analogue for surfaces of nonpositive of a theorem from~\cite{AN2}. As the geometric situation is slightly different from~\cite{AN2}, we will recall the main lines of the proof where the geometric properties appear and focus on the differences. We refer the reader to~\cite{AN2} for the details\footnote{We assume the reader is familiar with the proof of~\cite{AN2}.}. On~\cite{AN2}, the proof of the analogue of theorem~\ref{normestbis} (section $3$ and more precisely corollary $3.5$) relies on a study of the action of $P_{\alpha_{k}}U^{\eta}P_{\alpha_{k-1}}\cdots U^{\eta}P_{\alpha_0}$ on a particular family of Lagrangian states. This reduction was possible because of the introduction of the cutoffs operators $\Op_{\hbar}(\chi^{(k)})$ (see section $3$ in~\cite{AN2} for the details).

\subsection{Evolution of a WKB state}

Consider $u_{\hbar}(0,x)=a_{\hbar}(0,x)e^{\frac{\imath}{\hbar}S(0,x)}$ a Lagrangian state, where $a_{\hbar}(0,\bullet)$ and $S(0,\bullet)$ are smooth functions on a subset $\Omega$ in $M$ and $a_{\hbar}(0,\bullet)\sim\sum_k\hbar^ka_k(0,\bullet)$. This represents a Lagrangian state which is supported on the Lagrangian manifold $\mathcal{L}(0):=\{(x,d_xS(0,x):x\in\Omega\}$. According to~\cite{AN2}, if we are able to understand the action of $P_{\alpha_{k}}U^{\eta}P_{\alpha_{k-1}}\cdots U^{\eta}P_{\alpha_0}$ on Lagrangian states (with specific initial Lagrangian manifolds: see next paragraph), then we can derive our main theorem. A strategy to estimate this action is to use a WKB Ansatz. Recall that if we note $\tilde{u}(t):=U^tU(0)$, then, for any integer $N$, the state $\tilde{u}(t)$ can be approximated to order $N$ by a Lagrangian state $u(t)$ of the form
$$u(t,x):=e^{\frac{\imath}{\hbar}S(t,x)}a_{\hbar}(t,x)=e^{\frac{\imath}{\hbar}S(t,x)}\sum_{K=0}^{N-1}\hbar^ka_k(t,x).$$
As $u$ is supposed to solve $\imath\hbar\frac{\Delta}{2} u=\partial_t u$ (up to an error term of order $N$), we know that $S(t,x)$ and the $a_k(t,x)$ satisfy several partial differential equations. In particular, $S(t,x)$ must solve the Hamilton-Jacobi equation
$$\frac{\partial S}{\partial t}+H(x,d_xS)=0.$$
Assume that, on a certain time interval (for instance $s\in[0, \eta]$), the above equations have a
well defined smooth solution $S(s, x)$, meaning that the transported Lagrangian manifold
$\mathcal{L}(s) = g^s\mathcal{L}(0)$ is of the form $\mathcal{L}(s) = \{(x, d_xS(s, x))\}$, where $S(s)$ is a smooth function
on the open set $\pi\mathcal{L}(s)$.\\
As in~\cite{AN2}, we shall say that a Lagrangian manifold $\mathcal{L}$ is "projectible" if the
projection $\pi: \mathcal{L}\rightarrow M$ is a diffeomorphism onto its image. If the projection of $\mathcal{L}$ to $M$
is simply connected, this implies $\mathcal{L}$ is the graph of $dS$ for some function $S$: we say that $\mathcal{L}$ is generated by $S$.\\
Suppose now that, for $s\in[0,\eta]$, the Lagrangian $\mathcal{L}(s)$ is "projectible". Then, this family of Lagrangian manifolds define an induced flow on $M$, i.e.
$$g^t_{S(s)}:x\in\pi\mathcal{L}(s)\mapsto\pi g^t (x,d_xS(s,x))\in\pi\mathcal{L}(s+t).$$
This flow satifies a property of semi-group as follows: $g^t_{S(s+\tau)}\circ g_{S(s)}^{\tau}=g_{S(s)}^{t+\tau}$. Using this flow, we define an operator that sends functions on $\pi\mathcal{L}(s)$ into functions on $\pi\mathcal{L}(s+t)$:
$$T^t_{S(s)}(a)(x):=a\circ g_{S(s+t)}^{-t}(x)\left(J^{-t}_{S(s+t)}(x)\right)^{\frac{1}{2}},$$
where $J^{t}_{S(s)}(x)$ is the Jacobian of the map $g^t_{S(s)}$ at point $x$ (w.r.t. the riemannian volume). This operator allows to give an explicit expression for all the $a_k(t)$~\cite{AN2}, i.e.
$$a_k(t):=T^t_{S(0)}a_0(0)\ \text{and}\ a_k(t):=T^t_{S(0)}a_k(0)+\int_0^tT^{t-s}_{S(s)}\left(\frac{\imath\Delta a_{k-1}(s)}{2}\right)ds.$$
Regarding the details of the proof in~\cite{AN2}, we know that there are two main points where the dynamical properties of the manifold are used: 
\begin{itemize}
\item the evolution of the Lagrangian manifold under the action of $P_{\alpha_{k}}U^{\eta}P_{\alpha_{k-1}}\cdots U^{\eta}P_{\alpha_0}$ (section $3.4.1$ in~\cite{AN2});
\item the value of $J^{t}_{S(0)}$ for large $t$ (section $3.4.2$ in~\cite{AN2}).
\end{itemize}
We will discuss these two points in the two following paragraphs. We will recall what was proved for these two questions in section $3.4$ of~\cite{AN2} and see how it can be translated in the setting of surfaces of nonpositive curvature.

\subsection{Evolution of the Lagrangian manifolds}

The first thing we need to understand is how the Lagrangian manifolds evolve under the action of the operator $P_{\alpha_{k}}U^{\eta}P_{\alpha_{k-1}}\cdots U^{\eta}P_{\alpha_0}$. According to~\cite{AN2}, we know that the introduction of the cutoff operator $\Op_{\hbar}(\chi)$ implies that we can restrict our selves to a particular family of Lagrangian states. Precisely, we fix some small parameter $\eta_1$ and we know that they must be localized on a piece of Lagrangian manifold $\mathcal{L}^0(0)$ which is included in the set $\cup_{|\tau|\leq\eta}g^{\tau}S^*_{z,\eta_1}M$ (where $S^*_{z,\eta_1}M:=\{(z,\xi):\|\xi\|_z^2=1+2\eta_1\}$). If we follow the method developped in~\cite{AN2}, we are given a sequence of Lagrangian manifolds $\mathcal{L}^j(0)$ as follows:
$$\forall t\in[0,\eta],\ \forall j,\ \mathcal{L}^0(t):=g^t\mathcal{L}^0(0)\ \text{and}\ \mathcal{L}^j(t):=g^t\left(\mathcal{L}^{j-1}(\eta)\cap T^*\Omega_{\alpha_j}\right).$$
The manifold $\mathcal{L}^j(0)$ is obtained after performing $P_{\alpha_{j}}U^{\eta}P_{\alpha_{k-1}}\cdots U^{\eta}P_{\alpha_0}$ on the initial Lagrangian state. To show that the procedure from~\cite{AN2} is consistent (i.e. performing several WKB Ansatz), we need to verify that the Lagrangian manifold $\mathcal{L}^j(t)$ does not develop caustics and remains "projectible". The only geometric properties which were used to derive these two properties were:
\begin{itemize}
\item $M$ has no conjugate points (to derive that $S^j$ will not develop caustics); 
\item the injectivity radius is larger than $2$ (to ensure the "projectible" property).
\end{itemize}
In our setting, these two properties remain true (in particular, a surface of nonpositive curvature has no conjugate points~\cite{Ru}). Finally, we undeline that, thanks to the construction of the strong unstable foliation for surfaces of nonpositive curvature, any vector in $S^*_{z,\eta_1}M$ becomes uniformly close to the unstable subspace under the action of $d_{\rho}g^t$ (see lemma~\ref{unstable}). As a consequence, under the geodesic flow, a piece of sphere becomes uniformly close to the unstable foliation as $j$ tends to infinity. This point is the main difference with~\cite{AN2}. In fact, if we consider an Anosov geodesic flow, we have the stronger property that a piece of sphere becomes exponentially close to the unstable foliation, as $j$ tends to infinity. However, we will check that this property is sufficient for our needs.
\begin{rema} At this point of the proof, we can ask about an extension of these results to manifolds without conjugate points. According to~\cite{Ru}, the `uniform divergence' property (given by lemma~\ref{unstable}) is true for surfaces without conjugate points. We mention that this property fails in higher dimension for manifolds without conjugate points.
\end{rema}

\subsection{Estimates on the induced Jacobian}

As was already mentioned, the Jacobian $J^t_{S^j}$ of the map $g^t_{s_j}$ appears in the WKB expansion of a Lagrangian state evolved under the operator $P_{\alpha_{j}}U^{\eta}P_{\alpha_{j-1}}\cdots U^{\eta}P_{\alpha_0}$. Precisely, by iterating the WKB Ansatz, we have to estimate the following quantity (see equation $3.22$ in~\cite{AN2}):
\begin{equation}\label{coarsejac}J_k(x):=\left(J^{-\eta}_{S^{k-1}}(x)J_{S^{k-2}}^{-\eta}(g^{-\eta}_{S^k}(x))\cdots J_{S^1}^{-\eta}(g_{S^k}^{(-k+2)\eta}(x))\right)^{\frac{1}{2}}.\end{equation}
This Jacobian appears in each term of the WKB expansion of a Lagrangian state evolved under the operator $P_{\alpha_{k}}U^{\eta}P_{\alpha_{k-1}}\cdots U^{\eta}P_{\alpha_0}$ (see the formulas for the $a_p$). It is necessary to provide a way to bound this quantity as it will appear in the control of every derivatives of the WKB expansion. According to the proof in~\cite{AN2}, if we are able to bound uniformly this quantity, the bound we will obtain is the one that will appear in theorem~\ref{normestbis}. This point of the proof is the main difference with the proof in the Anosov case. So, our goal in this paragraph is to provide an upper bound on~(\ref{coarsejac}). This last quantity can be rewritten
$$J_k(x):=\exp\left(\frac{1}{2}\left(\log J^{-\eta}_{S^{k-1}}(x)+\log J_{S^{k-2}}^{-\eta}(g^{-\eta}_{S^k}(x))\cdots +\log J_{S^1}^{-\eta}(g_{S^k}^{(-k+2)\eta}(x))\right)\right).$$
As the Lagrangian $\mathcal{L}^j$ become uniformly close to the unstable foliation when $j$ tends to infinity, we know that, for every $\varepsilon'>0$, there exists some integer $j(\eta,\varepsilon')$ such that
$$\forall j\geq j(\eta,\varepsilon'),\ \forall \rho=(x,\xi)\in\mathcal{L}^j(0),\ |\log J_{S^{j}}^{-\eta}(x)-\log J_{S^u(\rho)}^{-\eta}(x)|\leq \varepsilon',$$
where $S^u(\rho)$ generates the local unstable manifold at point $\rho$ (which is a Lagrangian submanifold). Therefore, we find that there exists a constant $C(\varepsilon',\eta)$ (depending only on $\varepsilon'$ and $\eta$) such that, uniformly with respect to $k$ and to $\rho$ in $\mathcal{L}^k(0)$,
$$J_k(x)\leq C(\varepsilon',\eta)e^{k\varepsilon'}\prod_{j=0}^{k-1}J^{-\eta}_{S^u(g^{(-j+1)\eta}\rho)}(g_{S^k}^{(-j+1)\eta}(x))=C(\varepsilon',\eta)e^{k\varepsilon'}J^{(1-k)\eta}_{S^u(\rho)}(x).$$
The Jacobian $J^{-\eta}_{S^u(\rho)}$ measures the contraction of $g^{-\eta}$ along the unstable direction. From the construction of the unstable Riccati solution $U^u_{\rho}(s)$, we know that $U^u_{\rho}(s)$ also measures the contraction of $g^{-\eta}$ along $E^u(\rho)$. In fact, according to section~\ref{nonpositive}, one has
$$\|d_{\rho}g^{-t}_{|E^u(\rho)}\|\leq\sqrt{1+K_0}e^{\int_0^{-t}U^u_{\rho}(s)ds}.$$
As a consequence, there exists an uniform constant $C$ (depending only on the manifold) such that:
$$J^{(1-k)\eta}_{S^u(\rho)}(x)\leq C e^{\int_0^{(1-k)\eta}U^u_{\rho}(s)ds}.$$
Using then relation~(\ref{continuity2}) between the discrete Riccati solution $f$ and the continuous one, we find that there exists a constant $C(\epsilon,\eta,\epsilon_0)$ such that, uniformly in $k$,
$$\sup_{x\in\pi\mathcal{L}^k(0)}J_k(x)\leq C(\epsilon,\eta,\epsilon_0)e^{2k\eta\epsilon\epsilon_0}\exp\left(-\frac{1}{2}\sum_{j=0}^{k-1}f(\sigma^j\alpha)\right).$$
Finally, this last inequality gives us a bound on the quantity~(\ref{coarsejac}). This estimate is not as sharp as the one derived in~\cite{AN2} (equation $3.23$ for instance) however it is sufficient as the correction term is not too large: it is of order $\hbar^{-\epsilon}$.
\begin{rema}
We underline that we used the continuity of $U^u$ to go from the continuous representation of the upper bound of $J_k$ to the one in terms of the discrete Riccati solution. We underline again that this property fails for surfaces without conjugate points~\cite{BBB}.
\end{rema}

\section{Applying the uncertainty principle for quantum pressures}\label{s:appl-uncert}

In this section, we would like to prove inequality~(\ref{mainest-nonpos}) which was a crucial step of our proof. To do this, we follow the same lines as in~\cite{GR} (section $5.3$) and prove the following proposition:
\begin{prop}\label{lowerbound} With the notations of section~\ref{s:outline}, one has:
\begin{equation}\label{estQentr}
p\left(\overline{\mu}_{\hbar}^{\overline{\Sigma}_+},\overline{\mathcal{C}}_{\hbar}^+\right)+p\left(\overline{\mu}_{\hbar}^{\overline{\Sigma}_-},\overline{\mathcal{C}}^-_{\hbar}\right)\geq -\log C-(1+\epsilon'+4\epsilon)n_E(\hbar),\end{equation}
where $p$ is defined by~(\ref{pressure}) and where $C\in\mathbb{R}_+^*$ does not depend on $\hbar$ (but depends on the other parameters $(\epsilon,\epsilon_0,\eta)$).
\end{prop}
To prove this result, we will proceed in three steps. First, we will introduce an energy cutoff in order to get the sharpest bound as possible in our application of the uncertainty principle. Then, we will apply the uncertainty principle and derive a lower bound on $p\left(\overline{\mu}_{\hbar}^{\overline{\Sigma}_+},\overline{\mathcal{C}}_{\hbar}^+\right)+p\left(\overline{\mu}_{\hbar}^{\overline{\Sigma}_-},\overline{\mathcal{C}}^-_{\hbar}\right)$. Finally, we will use sharp estimates of theorem~\ref{normestbis} to conclude.

\subsection{Energy cutoff}

\label{cut}

Before applying the uncertainty principle, we proceed to sharp energy cutoffs so as to get precise lower bounds on the quantum pressure (as it was done in \cite{An}, \cite{AN2} and \cite{AKN}). These cutoffs are made in our microlocal analysis in order to get as good exponential decrease as possible of the norm of the refined quantum partition. This cutoff in energy is possible because even if the distributions $\mu_{\hbar}$ are defined on $T^*M$, they concentrate on the energy layer $S^*M$. The following energy localization is made in a way to compactify the phase space and in order to preserve the semiclassical measure.\\
Let $\delta_0$ be a positive number less than $1$ and $\chi_{\delta_0}(t)$ in $\mathcal{C}^{\infty}(\mathbb{R},[0,1])$. Moreover, $\chi_{\delta_0}(t)=1$ for $|t|\leq e^{-\delta_0 /2}$ and $\chi_{\delta_0}(t)=0$ for $|t|\geq 1$. As in~\cite{AN2}, the sharp $\hbar$-dependent cutoffs are then defined in the following way:
$$\forall\hbar\in(0,1),\ \forall n\in\mathbb{N},\ \forall\rho\in T^*M,\ \ \ \ \chi^{(n)}(\rho,\hbar):=\chi_{\delta_0}(e^{-n\delta_0}\hbar^{-1+\delta_0}(H(\rho)-1/2)).$$
For $n$ fixed, the cutoff $\chi^{(n)}$ is localized in an energy interval of length $2e^{n\delta_0}\hbar^{1-\delta_0}$ centered around the energy layer $\mathcal{E}$. In this paper, indices $n$ will satisfy $2e^{n\delta_0}\hbar^{1-\delta_0}<<1$. It implies that the widest cutoff is supported in an energy interval of microscopic length and that $n\leq K_{\delta_0}|\log\hbar|$, where $K_{\delta_0}\leq\delta_0^{-1}$. Using then a non standard pseudodifferential calculus (see~\cite{AN2} for a brief reminder of the procedure from~\cite{SZ}), one can quantize these cutoffs into pseudodifferential operators. We will denote $\Op(\chi^{(n)})$ the quantization of $\chi^{(n)}$. The main properties of this quantization are recalled in the appendix of~\cite{GR}. In particular, the quantization of these cutoffs preserves the eigenfunctions of the Laplacian:
\begin{prop}\label{hyp}~\cite{AN2} For any fixed $L>0$, there exists $\hbar_L$ such that for any $\hbar\leq\hbar_L$, any $n\leq K_{\delta}|\log\hbar|$ and any sequence $\beta$ of length $n$, the Laplacian eigenstate verify
$$\left\|\left(1-\Op\left(\chi^{(n)}\right)\right)\pi_{\beta}\psi_{\hbar}\right\|\leq \hbar^L\|\psi_{\hbar}\|.$$
\end{prop}

\subsection{Applying theorem~\ref{t:uncert}}
\label{several}
 Let $\|\psi_{\hbar}\|=1$ be a fixed element of the sequence of eigenfunctions of the Laplacian defined earlier, associated to the eigenvalue $-\frac{1}{\hbar^2}$.\\
To get bound on the pressure of the suspension measure, the uncertainty principle should not be applied to the eigenvectors $\psi_{\hbar}$ directly but it will be applied several times. Precisely, we will apply it to each $P_{\gamma}\psi_{\hbar}:=P_{\gamma_1}P_{\gamma_0}(-\eta)\psi_{\hbar}$ where $\gamma=(\gamma_0,\gamma_1)$ varies in $\{1,\cdots, K\}^2$. In order to apply the uncertainty principle to $P_{\gamma}\psi_{\hbar}$, we introduce new families of quantum partitions corresponding to each $\gamma$.\\
Let $\gamma=(\gamma_0,\gamma_1)$ be an element of $\{1,\cdots,K\}^2$. We define $\gamma.\alpha'=(\gamma_0,\gamma_1,\alpha')$. Introduce the following families of indices:
$$I_{\hbar}(\gamma):=\left\{(\alpha'):\gamma.\alpha'\in I^{\eta}(\hbar)\right\},$$
$$K_{\hbar}(\gamma):=\left\{(\beta'):\beta'.\gamma\in K^{\eta}(\hbar)\right\}.$$
We underline that each sequence $\alpha$ of $I^{\eta}(\hbar)$ can be written under the form $\gamma.\alpha'$ where $\alpha'\in I_{\hbar}(\gamma)$. The same works for $K^{\eta}(\hbar)$. The following partitions of identity can be associated to these new families, for $\alpha'\in I_{\hbar}(\gamma)$ and $\beta'\in K_{\hbar}(\gamma)$,
$$\tilde{\tau}_{\alpha'}=P_{\alpha_n'}(n\eta)\cdots P_{\alpha_2'}(2\eta),$$
$$\tilde{\pi}_{\beta'}=P_{\beta_{-n}'}(-n\eta)\cdots P_{\beta_{-2}'}(-2\eta).$$
The families $(\tilde{\tau}_{\alpha'})_{\alpha'\in I_{\hbar}(\gamma)}$ and $(\tilde{\pi}_{\beta'})_{\beta'\in I_{\hbar}(\gamma)}$ form quantum partitions of identity~\cite{GR}.\\
Given these new quantum partitions of identity, the unceratinty principle should be applied for given initial conditions $\gamma=(\gamma_0,\gamma_1)$ in times $0$ and $1$. We underline that for $\alpha'\in I_{\hbar}(\gamma)$ and $\beta'\in K_{\hbar}(\gamma)$:
\begin{equation}\label{switch}\tilde{\tau}_{\alpha'}U^{-\eta}P_{\gamma}=\tau_{\gamma.\alpha'}U^{-\eta}\ \text{and}\ \tilde{\pi}_{\beta'}P_{\gamma}=\pi_{\beta'.\gamma},\end{equation}
where $\gamma.\alpha'\in I^{\eta}(\hbar)$ and $\beta'.\gamma\in K^{\eta}(\hbar)$ by definition. In equality~(\ref{switch}) appears the fact that the definitions of $\tau$ and $\pi$ are slightly different (see~(\ref{taubis}) and~(\ref{pibis})). It is due to the fact that we want to compose $\tilde{\tau}$ and $\tilde{\pi}$ with the same operator $P_{\gamma}$.\\
Suppose now that $\|P_{\gamma}\psi_{\hbar}\|$ is not equal to $0$. We apply the quantum uncertainty principle~\ref{t:uncert} using that
\begin{itemize}
\item $(\tilde{\tau}_{\alpha'})_{\alpha'\in I_{\hbar}(\gamma)}$ and $(\tilde{\pi}_{\beta'})_{\beta'\in K_{\hbar}(\gamma)}$ are partitions of identity;
\item the cardinal of $I_{\hbar}(\gamma)$ and $K_{\hbar}(\gamma)$ is bounded by $\mathcal{N}\simeq \hbar^{-K_0}$ where $K_0$ is some fixed positive number (depending on the cardinality of the partition $K$, on $a_0$, on $b_0$ and $\eta$);
\item $\Op(\chi^{(k')})$ is a family of bounded bounded operators $O_{\beta'}$ (where $k'$ is the length of $\beta'$);
\item the constants $W_{\gamma.\alpha'}^+$ and $W_{\beta.\gamma}^-$ are bounded by $\hbar^{-\frac{b_0}{2\epsilon_0}}$;
\item the parameter $\delta'$ can be taken equal to $\|P_{\gamma}\psi_{\hbar}\|^{-1}\hbar^L$ where $L$ is such that $\hbar^{L-K_0-\frac{b_0}{2\epsilon_0}}\ll e^{2k\eta\epsilon\epsilon_0} \hbar^{-1/2-c\delta_0}$ for every $k\ll \frac{1}{\epsilon\eta}|\log\hbar|$ (see proposition~\ref{hyp} and the upper bound in theorem~\ref{normestbis});
\item $U^{-\eta}$ is an isometry;
\item $\tilde{\psi_{\hbar}}:=\frac{P_{\gamma}\psi_{\hbar}}{\|P_{\gamma}\psi_{\hbar}\|}$ is a normalized vector. 
\end{itemize}
Applying the uncertainty principle~\ref{t:uncert} for quantum pressures, one gets:
\begin{coro}\label{UP} Suppose that $\|P_{\gamma}\psi_{\hbar}\|$ is not equal to $0$. Then, one has
$$p_{\tilde{\tau}}(U^{-\eta}\tilde{\psi}_{\hbar})+p_{\tilde{\pi}}(\tilde{\psi}_{\hbar})\geq-2\log \left(c_{\chi}^{\gamma}(U^{-\eta})+\hbar^{L-K_0-\frac{b_0}{2\epsilon_0}}\|P_{\gamma}\psi_{\hbar}\|^{-1}\right),$$
where  $\displaystyle c_{\chi}^{\gamma}(U^{-\eta})=\max_{\alpha'\in I_{\hbar}(\gamma),\beta'\in K_{\hbar}(\gamma)}\left(W_{\gamma.\alpha'}^+W_{\beta'.\gamma}^-\|\tilde{\tau}_{\alpha'}U^{-\eta}\tilde{\pi}_{\beta'}^* \Op(\chi^{(k')})\|\right).$
\end{coro}
Under this form, the quantity $\|P_{\gamma}\psi_{\hbar}\|^{-1}$ appears several times and we would like to get rid of it. First, remark that the quantity $c_{\chi}^{\gamma}(U^{-\eta})$ can be easily replaced by \begin{equation}\label{norm}c_{\chi}(U^{-\eta}):=\max_{\gamma\in\{1,\cdots,K\}^2}\max_{\alpha'\in I_{\hbar}(\gamma),\beta'\in K_{\hbar}(\gamma)}\left(W_{\gamma.\alpha'}^+W_{\beta'.\gamma}^-\|\tilde{\tau}_{\alpha'}U^{-\eta}\tilde{\pi}_{\beta'}^* \Op(\chi^{(k')})\|\right),\end{equation}
which is independent of $\gamma$. Then, one has the following lower bound:
\begin{equation}\label{bound}-2\log \left(c_{\chi}^{\gamma}(U^{-\eta})+\hbar^{L-K_0}\|P_{\gamma}\psi_{\hbar}\|^{-1}\right)\geq-2\log \left(c_{\chi}(U^{-\eta})+\hbar^{L-K_0-\frac{b_0}{2\epsilon_0}}\right)+2\log\|P_{\gamma}\psi_{\hbar}\|^{2}.\end{equation}
as $\|P_{\gamma}\psi_{\hbar}\|\leq1$. Now that we have given an alternative lower bound, we rewrite the entropy term $h_{\tilde{\tau}}(U^{-\eta}\tilde{\psi}_{\hbar})$ of the quantum pressure $p_{\tilde{\tau}}(U^{-\eta}\tilde{\psi}_{\hbar})$ as follows:
$$h_{\tilde{\tau}}(U^{-\eta}\tilde{\psi}_{\hbar})=-\sum_{\alpha'\in I_{\hbar}(\gamma)}\|\tilde{\tau}_{\alpha'}U^{-\eta}\tilde{\psi}_{\hbar}\|^2\log\|\tilde{\tau}_{\alpha'}U^{-\eta}P_{\gamma}\psi_{\hbar}\|^2+\sum_{\alpha'\in I_{\hbar}(\gamma)}\|\tilde{\tau}_{\alpha'}U^{-\eta}\tilde{\psi}_{\hbar}\|^2\log\|P_{\gamma}\psi_{\hbar}\|^2.$$
Using the fact that $\psi_{\hbar}$ is an eigenvector of $U^{\eta}$ and that $(\tilde{\tau}_{\alpha'})_{\alpha'\in I_{\hbar}(\gamma)}$ is a partition of identity, one has:
$$h_{\tilde{\tau}}(U^{-\eta}\tilde{\psi}_{\hbar})=-\frac{1}{\|P_{\gamma}\psi_{\hbar}\|^2}\sum_{\alpha'\in I_{\hbar}(\gamma)}\|\tau_{\gamma.\alpha'}\psi_{\hbar}\|^2\log\|\tau_{\gamma.\alpha'}\psi_{\hbar}\|^2+\log\|P_{\gamma}\psi_{\hbar}\|^2.$$
The same holds for the entropy term $h_{\tilde{\pi}}(\tilde{\psi}_{\hbar})$ of the quantum pressure $p_{\tilde{\pi}}(\tilde{\psi}_{\hbar})$ (using here equality~(\ref{switch})):
$$h_{\tilde{\pi}}(\tilde{\psi}_{\hbar})=-\frac{1}{\|P_{\gamma}\psi_{\hbar}\|^2}\sum_{\beta'\in K_{\hbar}(\gamma)}\|\pi_{\beta'.\gamma}\psi_{\hbar}\|^2\log\|\pi_{\beta'.\gamma}\psi_{\hbar}\|^2+\log\|P_{\gamma}\psi_{\hbar}\|^2.$$
Combining these last two equalities with~(\ref{bound}), we find that
\begin{equation}\label{inter}
-\sum_{\alpha'\in I_{\hbar}(\gamma)}\|\tau_{\gamma.\alpha'}\psi_{\hbar}\|^2\log\|\tau_{\gamma.\alpha'}\psi_{\hbar}\|^2-2\sum_{\alpha'\in I_{\hbar}(\gamma)}\|\tau_{\gamma.\alpha'}\psi_{\hbar}\|^2\log W_{\gamma.\alpha'}^+\end{equation}
$$-\sum_{\beta'\in K_{\hbar}(\gamma)}\|\pi_{\beta'.\gamma}\psi_{\hbar}\|^2\log\|\pi_{\beta'.\gamma}\psi_{\hbar}\|^2-2\sum_{\beta'\in K_{\hbar}(\gamma)}\|\pi_{\beta'.\gamma}\psi_{\hbar}\|^2\log W^-_{\beta'.\gamma}\geq-2\|P_{\gamma}\psi_{\hbar}\|^{2}\log \left(c_{\chi}(U^{-\eta})+\hbar^{L-K_0-\frac{b_0}{2\epsilon_0}}\right).$$
This expression is very similar to the definition of the quantum pressure. We also underline that this lower bound is trivial in the case where $\|P_{\gamma}\psi_{\hbar}\|$ is equal to $0$. Using the following numbers:
\begin{equation}\label{cgamma}c_{\gamma.\alpha'}=c_{\beta'.\gamma}=c_{\gamma}=\frac{f(\gamma)}{\sum_{\gamma'\in\{1,\cdots,K\}^2}f(\gamma')\|P_{\gamma'}\psi_{\hbar}\|^2},\end{equation}
one can derive, as in~\cite{GR}, the following property:
\begin{coro}\label{corouncert}
One has:
\begin{equation}\label{uncerteigenfct}p\left(\overline{\mu}_{\hbar}^{\overline{\Sigma}_+},\overline{\mathcal{C}}_{\hbar}^+\right)+p\left(\overline{\mu}_{\hbar}^{\overline{\Sigma}_-},\overline{\mathcal{C}}^-_{\hbar}\right)\geq-2\log \left( c_{\chi}(U^{-\eta})+\hbar^{L-K_0-\frac{b_0}{2\epsilon_0}}\right)-\log\left(\max_{\gamma}c_{\gamma}\right).\end{equation}
\end{coro}
As expected, by a careful use of the entropic uncertainty principle, we have been able to obtain a lower bound on the pressures of the measures $\overline{\mu}_{\hbar}^{\overline{\Sigma}_+}$ and $\overline{\mu}_{\hbar}^{\overline{\Sigma}_-}$.

\subsection{The conclusion} To conclude the proof of proposition~\ref{lowerbound}, we use theorem~\ref{normestbis} to give an upper bound on $c_{\chi}(U^{-\eta})$. From our assumption on $L$, we know that $\hbar^{L-K_0-\frac{b_0}{2\epsilon_0}}\ll c_{\chi}(U^{-\eta}).$ As $k\eta\leq n_E(\hbar)/\epsilon_0$, we also have that
$$c_{\chi}(U^{-\eta})\leq C_{\mathcal{K}}(\epsilon,\eta,\epsilon_0)\hbar^{-\frac{1}{2}-c\delta_0}e^{4\epsilon n_E(\hbar)}.$$
For $\delta_0$ small enough, we find the expected property.$\square$

\appendix

\section{Uncertainty principle for the quantum pressure}
\label{quant-pressure}
In~\cite{AN2}, generalizations of the entropic uncertainty principle were derived for quantum pressures. We saw that the use of this thermodynamic formalism was crucial in our proof and we recall in this section the main results from~\cite{AN2} (section $6$) on quantum pressures. Consider two partitions of identity $(\pi_k)_{k=1}^{\mathcal{N}}$ and $(\tau_j)_{j=1}^{\mathcal{M}}$ on $L^2(M)$, i.e.
$$\sum_{k=1}^{\mathcal{N}}\pi_k^*\pi_k=\text{Id}_{L^2(M)}\ \text{and}\ \sum_{j=1}^{\mathcal{M}}\tau_j^*\tau_j=\text{Id}_{L^2(M)}.$$
We also introduce two families of positive numbers: $(V_k)_{k=1}^{\mathcal{N}}$ and $(W_j)_{j=1}^{\mathcal{M}}$. We denote $A:=\max_kV_k$ and $B:=\max_jW_j.$ One can then introduce the quantum pressures associated to these families, for a normalized vector $\psi$ in $L^2(M)$,
$$p_{\pi}(\psi):=-\sum_{k=0}^{\mathcal{N}}\|\pi_k\psi\|^2_{L^2(M)}\log\|\pi_k\psi\|^2_{L^2(M)}-2\sum_{k=0}^{\mathcal{N}}\|\pi_k\psi\|^2_{L^2(M)}\log V_k$$
and
$$p_{\tau}(\psi):=-\sum_{j=0}^{\mathcal{M}}\|\tau_j\psi\|^2_{L^2(M)}\log\|\tau_j\psi\|^2_{L^2(M)}-2\sum_{j=0}^{\mathcal{M}}\|\tau_j\psi\|^2_{L^2(M)}\log W_j.$$
The main result on these quantities that was derived in~\cite{AN2} was theorem $6.5$:
\begin{theo}\label{t:uncert} Under the previous setting, suppose $\mathcal{U}$ is an isometry of $L^2(M)$ and suppose $(O_{k})_{k=1}^{\mathcal{N}}$ is a family of bounded operators. Let $\delta'$ be a positive number and $\psi$ be a vector in $\mathcal{H}$ of norm $1$ such that
$$\|(Id-O_{k})\pi_{k}\psi\|_{L^2(M)}\leq\delta'.$$
Then, one has
$$p_{\tau}(\mathcal{U}\psi)+p_{\pi}(\psi)\geq-2\log \left(c_O^{\alpha,\beta}(\mathcal{U})+\mathcal{N}AB\delta'\right),$$
where $c_O^{\alpha,\beta}(\mathcal{U}):=\sup_{j,k}\{V_kW_j\|\tau_j\mathcal{U}\pi_k^*O_k\|\}.$
\end{theo}

\section{Kolmogorov-Sinai entropy}
\label{KSentropy}
Let us recall a few facts about Kolmogorov-Sinai (or metric) entropy that can be found for example in~\cite{Wa}. Let $(X,\mathcal{B},\mu)$ be a measurable probability space, $I$ a finite set and $P:=(P_{\alpha})_{\alpha\in I}$ a finite measurable partition of $X$, i.e. a finite collection of measurable subsets that forms a partition. Each $P_{\alpha}$ is called an atom of the partition. Assuming $0\log 0=0$, one defines the entropy of the partition as:
\begin{equation}\label{defent0}H(\mu,P):=-\sum_{\alpha\in I}\mu(P_{\alpha})\log\mu(P_{\alpha})\geq 0.\end{equation}
Given two measurable partitions $P:=(P_{\alpha})_{\alpha\in I}$ and $Q:=(Q_{\beta})_{\beta\in K}$, one says that $P$ is a refinement of $Q$ if every element of $Q$ can be written as the union of elements of $P$ and it can be shown that $H(\mu,Q)\leq H(\mu,P)$. Otherwise, one denotes $P\vee Q:=(P_{\alpha}\cap Q_{\beta})_{\alpha\in I,\beta\in K}$ their join (which is still a partition) and one has $H(\mu,P\vee Q)\leq H(\mu,P)+H(\mu,Q)$ (subadditivity property). Let $T$ be a measure preserving transformation of $X$. The $n$-refined partition $\vee_{i=0}^{n-1}T^{-i}P$ of $P$ with respect to $T$ is then the partition made of the atoms $(P_{\alpha_0}\cap\cdots\cap T^{-(n-1)}P_{\alpha_{n-1}})_{\alpha\in I^n}$. We define the entropy with respect to this refined partition:
\begin{equation}\label{nKSentropy}H_{n}(\mu, T, P)=-\sum_{|\alpha|=n}\mu(P_{\alpha_0}\cap\cdots\cap T^{-(n-1)}P_{\alpha_{n-1}})\log\mu(P_{\alpha_0}\cap\cdots\cap T^{-(n-1)}P_{\alpha_{n-1}}).\end{equation}
Using the subadditivity property of entropy, we have for any integers $n$ and $m$:
\begin{equation}\label{classicalsubad}
H_{n+m}(\mu, T, P)\leq H_{n}(\mu, T, P)+H_{m}(T^n\sharp\mu, T, P)=H_{n}(\mu, T, P)+H_{m}(\mu, T, P).
\end{equation}
For the last equality, it is important to underline that we really use the $T$-invariance of the measure $\mu$. A classical argument for subadditive sequences allows us to define the following quantity:
\begin{equation}\label{defentropy}h_{KS}(\mu,T,P):=\lim_{n\rightarrow\infty}\frac{H_n\left(\mu,T,P\right)}{n}.\end{equation}
It is called the Kolmogorov Sinai entropy of $(T,\mu)$ with respect to the partition $P$. The Kolmogorov Sinai entropy $h_{KS}(\mu,T)$ of $(\mu,T)$ is then defined as the supremum of $h_{KS}(\mu,T,P)$ over all partitions $P$ of $X$.


\begin{thebibliography}{99}
\bibitem{Ab} L.M.~Abramov \emph{On the entropy of a flow}, Translations of AMS $\mathbf{49}$, 167-170 (1966)
\bibitem{An} N.~Anantharaman \emph{Entropy and the localization of eigenfunctions}, Ann. of Math. $\mathbf{168}$, 435-475 (2008)
\bibitem{AKN} N.~Anantharaman, H.~Koch, S.~Nonnenmacher \emph{Entropy of eigenfunctions}, arXiv:0704.1564, International Congress of Mathematical Physics
\bibitem{AN2} N.~Anantharaman, S. Nonnenmacher \emph{Half-delocalization of eigenfunctions for the Laplacian on an Anosov manifold}, Ann. Inst. Fourier $\mathbf{57}$, 2465-2523 (2007)
\bibitem{BBB} W.~Ballmann, M.~Brin, K.~Burns \emph{On surfaces with no conjugate points}, Jour. Diff. Geom. $\mathbf{25}$, 249-273 (1987)
\bibitem{BP} L.~Barreira, Y.~Pesin \emph{Lectures on Lyapunov exponents and smooth ergodic theory}, Proc. of Symposia in Pure Math. $\mathbf{69}$, 3-89 (2001)
\bibitem{Burq} N.~Burq \emph{Mesures semi-classiques et mesures de défaut (d'après P.~Gérard, L.~Tartar et al.)}, Astérisque $\mathbf{245}$, 167-196, Séminaire Bourbaki, (1996-1997)
\bibitem{CdV} Y.~Colin de Verdière \emph{Ergodicité et fonctions propres du Laplacien}, Comm. in Math. Phys. $\mathbf{102}$, 497-502 (1985)
\bibitem{DS} M.~Dimassi, J.~Sj\"ostrand \emph{Spectral Asymptotics in the Semiclassical Limit} Cambridge University Press (1999)
\bibitem{Do} H.~Donnelly \emph{Quantum unique ergodicity}, Proc. of Amer. Math. Soc. $\mathbf{131}$, 2945-2951 (2002)
\bibitem{Eb0} P.~Eberlein \emph{When is a geodesic flow of Anosov type I}, Jour. Diff. Geom. $\mathbf{8}$, 437-463 (1973)
\bibitem{Eb} P.~Eberlein \emph{Geodesic flows in manifolds of nonpositive curvature}, Proc. of Symposia in Pure Math. $\mathbf{69}$, 525-571 (2001)
\bibitem{FM} A.~Freire, R.~Ma\~n\'e \emph{On the entropy of the geodesic flow for manifolds without conjugate points}, Inv. Math. $\mathbf{69}$, 375-392 (1982)
\bibitem{Gr} L.~Green \emph{Geodesic instability}, Proc. of Amer. Math. Soc. $\mathbf{7}$, 438-448 (1956)
\bibitem{Has} A.~Hassell \emph{Ergodic billiards that are not quantum unique ergodic. With an appendix by A. Hassell and L. Hillairet}, to appear in Ann. of Math.
\bibitem{LY} F.~Ledrappier, L.-S. Young \emph{The metric entropy of diffeomorphisms I. Characterization of measures satisfying Pesin's entropy formula}, Ann. of Math. $\mathbf{122}$, 509-539 (1985)
\bibitem{MU} H.~Maassen, J.B.~Uffink \emph{Generalized entropic uncertainty relations}, Phys. Rev. Lett. $\mathbf{60}$, 1103-1106 (1988)
\bibitem{GR} G.~Rivière \emph{Entropy of semiclassical measures in dimension 2}, arXiv:0809.0230 (2008)
\bibitem{RS} Z.~Rudnick, P. Sarnak \emph{The behaviour of eigenstates of arithmetic hyperbolic manifolds}, Comm. in Math. Phys. $\mathbf{161}$, 195-213 (1994)
\bibitem{R} D.~Ruelle \emph{An inequality for the entropy of differentiable maps}, Bol. Soc. Bras. Mat. $\mathbf{9}$, 83-87 (1978)
\bibitem{Ru} R.~O.~Ruggiero \emph{Dynamics and global geometry of manifolds without conjugate points}, Ensaios Mate. $\mathbf{12}$, Soc. Bras. Mate. (2007)
\bibitem{Sc} A.~Shnirelman \emph{Ergodic properties of eigenfunctions}, Usp. Math. Nauk. $\mathbf{29}$, 181-182 (1974)
\bibitem{SZ} J.~Sj\"ostrand, M.~Zworski \emph{Asymptotic distribution of resonances for convex obstacles}, Acta Math. $\mathbf{183}$, 191-253 (1999)
\bibitem{Wa} P.~Walters \emph{An introduction to ergodic theory}, Springer-Verlag, Berlin, New York (1982)
\bibitem{Ze} S.~Zelditch \emph{Uniform distribution of the eigenfunctions on compact hyperbolic surfaces}, Duke Math. Jour. $\mathbf{55}$, 919-941 (1987)
\end{thebibliography}
\end{document}